%% file: surface_arxiv_v2.tex
\documentclass[a4paper,11pt]{amsart}

\input{header-surfaces.tex}
\input{header-tikz-surfaces.tex}

\begin{document}

\title{Rigid divisors on surfaces}
\author{Andreas Hochenegger}
\author{David Ploog}

\thispagestyle{empty}

\begin{abstract}
We study effective divisors $D$ on surfaces with
$H^0(\OO_D)=\kk$ and $H^1(\OO_D)=H^0(\OO_D(D))=0$. We give a numerical
criterion for such divisors, following a general investigation of
negativity, rigidity and connectivity properties. Examples include
exceptional loci of rational singularities, and spherelike divisors.
\end{abstract}

\begingroup\renewcommand\thefootnote{}%
\addtocounter{footnote}{-2}\endgroup
\maketitle

\setcounter{tocdepth}{1}
\tableofcontents

\addtocontents{toc}{\protect\setcounter{tocdepth}{-1}}     

\section*{Introduction} 
\addtocontents{toc}{\protect\setcounter{tocdepth}{1}}      

\noindent
An important facet of algebraic geometry is the rich theory of positive divisors on varieties. The starting point are ample divisors, a notion which can be characterised geometrically, cohomologically and numerically. Moreover, there are various generalisations, leading to a mature tenet of positivity; see \cite{Lazarsfeld}. We point out that the theory was started by Zariski's work on surfaces.

In this article, we conduct a systematic study of negativity properties. Here, we restrict to surfaces but we feel that a numerical approach going beyond intersection numbers should yield results in higher dimension. Some classical negativity notions for surfaces are negative definiteness (appearing in Artin's contraction criterion) and Ramanujan's 1-connectedness.

Specifically, we consider effective divisors $D$ on a smooth algebraic surface $X$ which are \emph{well-connected} ($H^0(\OO_D)=\kk$) and \emph{rigid} as sheaves on $X$ (\ie $H^1(\OO_D)=H^0(\OO_D(D))=0$). 
We then call $D$ a \emph{$(-n)$-divisor} as necessarily $D^2 = -n < 0$.
This turns out to be a good notion which, for example, has a nice numerical characterisation.

Such divisors always consist of negative, rational curves. We don't discuss here negative curves on a fixed surface; for this, and the Bounded Negativity Conjecture, see \cite{BHKKMSRS}.

We have two principal motivations for this work. Geometrically, rigid and negative definite divisors arise as exceptional loci of rational singularities. In particular, $(-1)$-divisors come from blowing up smooth points, and numerical fundamental cycles of ADE singularities are $(-2)$-divisors. More generally, those cycles of rational singularities are always $(-n)$-divisors, see \autoref{prop:Znum-curvelike}.

Homologically, the structure sheaf $\OO_D$ of a $(-2)$-divisor is a 2-spherelike object in the sense of \cite{HKP}, \ie $\Hom(\OO_D,\OO_D)=\Ext^2(\OO_D,\OO_D)=\kk$ and $\Ext^1(\OO_D,\OO_D)=0$. This was our initial motivation to study these divisors.
In our previous work, to such $D$ we associate a natural, maximal subcategory of $\Db(X)$ in which $\OO_D$ becomes a 2-Calabi--Yau object. 
As a homological detour, we compute these spherical subcategories and the related asphericities for some examples in \autoref{sec:spherelike-divisors}.

We start out, in \autoref{sec:divisor-properties}, with a systematic study of divisors centered around negativity, rigidity and connectivity. Here we list the non-standard notions in a very terse fashion; for more elaborate definitions and comments, see the referenced subsections. Let $D$ be an effective divisor; all $C_i$ occurring below are reduced, irreducible curves.

\introitem{Negativity}{\autoref{sec:negativity-properties}}
$D$ is \emph{negatively closed} if $A^2<0$ for any $0 \led A \leqd D$.
It is \emph{negative definite} if $kD$ is negatively closed for all $k\geq1$.
It is \emph{negatively filtered} if there is $D=C_1+\cdots+C_m$ with $C_i.(C_i+\cdots+C_m)<0$ for $i=1,\ldots,m$.

\introitem{Rigidity}{\autoref{sec:rigidity-properties}}
$D$ is \emph{rigid as a subscheme} if $H^0(\OO_D(D))=0$.
It is \emph{Jacobi rigid} if $H^1(\OO_D)=0$.
It is \emph{rigid} if $D$ is rigid as a subscheme and Jacobi rigid, \ie $\OO_D$ is an infinitesimally rigid sheaf on $X$.

\introitem{Connectivity}{\autoref{sec:connectivity-properties}}
$D$ is \emph{well-connected} if $H^0(\OO_D)=\kk$.
It is \emph{1-connected} if $A.(D-A) \geq 1$ for all $0 \led A \led D$.
It is \emph{1-decomposable} if there is $D=C_1+\cdots+C_m$ with $C_i.(C_{i+1}+\cdots+C_m)=1$ for $i=1,\ldots,m-1$.

\smallskip

Our results about these properties are most concisely summed up in the following schematic. The ornaments $c$ and $n$ indicate that a property is \emph{closed under subdivisors} or \emph{numerical}, respectively; see \autoref{def:divisor-properties}.
By \autoref{prop:birational-properties}, all these properties are birationally invariant.
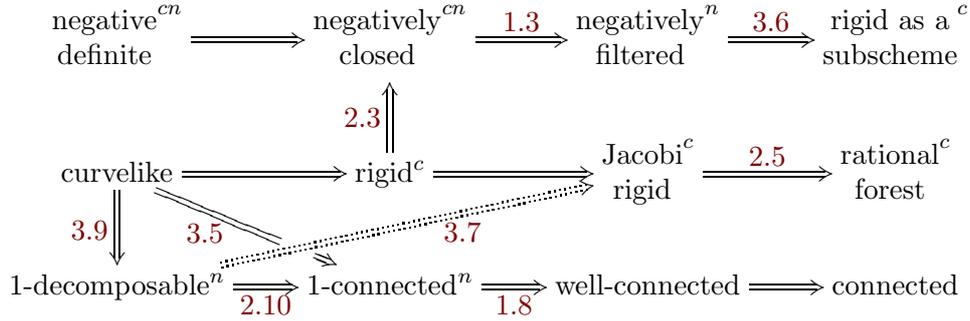
\[ 
\begin{tikzcd}
   \cdstext{negative}{definite}\closednum   \ar[r, Rightarrow] & 
   \cdstext{negatively}{closed}\closednum   \ar[r, Rightarrow, "\ref{lem:negatively-closed-implies-negatively-filtered}"] &
   \cdstext{negatively}{filtered}\num       \ar[r, Rightarrow, "\ref{lem:negatively-filtered-implies-subscheme-rigid}"] &
   \cdstext{rigid as a}{subscheme}\closed         
\\
   \cdstext{well-connected}{and rigid}                      \ar[dd, Rightarrow, "\ref{lem:has-1-decomposition}"']
                                        \ar[dr, Rightarrow, "\ref{lem:curvelike-num1conn}"'] \ar[r, Rightarrow ] & 
   \cdtext{rigid}\closed                   \ar[r, Rightarrow ] \ar[u, Rightarrow, "\ref{cor:rigid-totally}"] &
   \cdstext{Jacobi}{rigid}\closed           \ar[r, Rightarrow, "\ref{lem:rigid-tree-rat}"'] & 
   \cdstext{rational}{forest}\closed
\\
 & \cdtext{1-decomposable}\num             \ar[ru, dashed, Rightarrow, "\ref{lem:1-decomposition-given}"']
                                        \ar[d, dashed, Rightarrow, "\ref{lem:1-decomposition-given}"']
\\
   \cdtext{1-connected}\num                \ar[r, Rightarrow, "\ref{lem:1-connected-is-well-connected}"] &
   \cdtext{well-connected}            \ar[r, Rightarrow ] & 
   \cdtext{connected}  \\
&&&&
\end{tikzcd}
\]

The crucial new notion is that of 1-decomposability, because it enables a combinatorial grip on $H^0(\OO_D)=\kk$ and $H^1(\OO_D)=0$.
One of our main results is the following characterisation; see \autoref{thm:numerical-criterion}.

\begin{theoremalpha}
An effective divisor is well-connected and rigid if and only if it is rational, $1$-decomposable and negatively filtered.
\end{theoremalpha}

Let $D$ be a well-connected, rigid divisor, i.e.\ a $(-n)$-divisor with $n = -D^2$.
We are interested in simplifying $D$. Such a simplification can be
\begin{description}[leftmargin=2em]
\item[birational] if there is a contraction $\pi\colon X\to X'$ such that $D'=\pi(D)$ is a $(-n)$-divisor on $X'$ with $D = \pi^*D'$. For this to be possible, there has to be a $(-1)$-curve $E\subset X$ with $E.D=0$.
\item[homological] while $(-1)$-curves induce contractions, \ie categorical decompositions of the derived category, a $(-2)$-curve $C$ gives rise to an autoequivalence, the spherical twist $\TTT_{\OO_C}$. If $C\led D$ with $C.D=-1$, then $\TTT_{\OO_C(-1)}(\OO_D)=\OO_{D-C}$, and we say that $C$ \emph{can be twisted off} $D$. 
\item[numerical] a sum $D=D_1+\cdots+D_m$ with $(-n_i)$-divisors $D_i$ such that $D_i.(D_{i+1}+\cdots+D_m)=1$ and $-2 \geq -n_i \geq -n$ for all $i$ is called a \emph{curvelike decomposition}; see \autoref{def:curvelike-decomposition}.
\end{description}

There always is a decomposition in the numerical sense with strong properties; the following result is \autoref{prop:curvelike-chopped}.

\begin{theoremalpha}
\label{thm:B}
For $n\geq2$, any $(-n)$-divisor has a curvelike decomposition where all summands are pullbacks of either curves $C$ or chains of type $\intext{\tikzsimplechain{m}}$, such that $C^2,-m \in \{-2,\ldots,-n\}$.
\end{theoremalpha}

An alternative definition of $(-n)$-divisors $D$ is $H^\bullet(\OO_D)=H^\bullet(\OO_C)$ and $H^\bullet(\OO_D(D))=H^\bullet(\OO_C(C))$, where $C$ is a $(-n)$-curve, \ie a smooth rational curve such that $C^2=-n$, see \autoref{lem:deformations}. In \autoref{prop:Odual-hom-groups}, we show that likewise $H^\bullet(\OO_D(K_X))=H^\bullet(\OO_C(K_X))$ and $H^\bullet(\OO_D(D+K_X))=H^\bullet(\OO_C(C+K_X))$.

We say that $D$ is \emph{essentially a $(-n)$-curve}, if $D$ can be obtained from a $(-n)$-curve through blow-ups and spherical twists, see \autoref{def:minimal-divisor}. At the other extreme, we call $D$ a \emph{minimal $(-n)$-divisor} if neither operation is possible.
For general reasons, certain low-degree divisors are always essentially curves; see \autoref{cor:exceptional-divisors}, \autoref{cor:spherical-divisors-essential}, \autoref{prop:reduced-spherelike}, respectively.

\begin{theoremalpha}
A divisor which is exceptional, or spherical, or reduced and spherelike is essentially a $(-n)$-curve for some $n>0$.
\end{theoremalpha}

Moreover, the building blocks, \ie minimal $(-n)$-divisors, can be enumerated. For instance, there are five minimal $(-2)$-divisors on five curves; see \autoref{ex:minimally-spherelike-5}: 
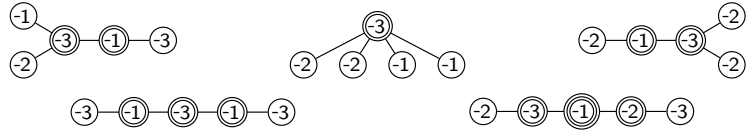
\begin{center}
  \hfill \tikztypeD{2}{1}{3,1,3}{2,2,1}
  \hfill \tikztreemult{2,2,1,1}{3}{1,1,1,1}{2}
  \hfill \tikztypeDmirror{2}{2}{3,1,2}{2,2,1}
  \hfill{ }

\medskip
  \hfill \tikzchainmult{3,1,3,1,3}{1,2,2,2,1}
  \hfill \tikzchainmult{2,3,1,2,3}{1,2,3,2,1}
  \hfill{ }
\end{center}

The question remains whether such graphs can be realised as dual intersection graphs of effective divisors. In \autoref{prop:existence-of-divisors}, we give a sufficient condition which is enough to deal with the above examples. This seems to be a subtle problem, and it would be interesting to study it in greater depth.

\subsubsection*{Conventions}
We work over an algebraically closed field $\kk$.
Throughout this article, we consider a smooth algebraic surface, usually denoted $X$, and effective projective divisors on it. For such divisors, the intersection product is well defined, see for example \cite[Ex.\ 2.4.9]{Fulton}. The canonical divisor on $X$ is written $K_X$ and, sometimes, just $K$. The structure sheaf of the surface is denoted $\OO$ instead of $\OO_X$.
Curves are always irreducible and reduced. By a $(-n)$-curve, we mean a smooth rational curve on $X$ of self-intersection number $-n$.

$\Hom$ and $\Ext^i$ refer to homorphism and extension spaces on the surface, unless we explicitly specify another variety.
Occasionally, we abbreviate dimensions $\hom(M,N) \coloneqq \dim\Hom(M,N)$ and $h^i(M) \coloneqq \dim H^i(M)$.
The $\Hom$ complex of two sheaves $M,N$ is $\Hom^\bullet(M,N) \coloneqq \bigoplus_i\Ext^i(M,N)[-i]$; it is a complex of vector spaces with zero differentials.

Sometimes we use $M' \into M \onto M''$ as a short hand for a short exact sequence.
Distinguished triangles in $\Db(X)$ are shown as $M'\to M\to M''$, omitting the connecting morphism $M''\to M'[1]$, and called just `triangles'.
We do not adorn the symbol for derived functors, \eg $f_*\colon \Db(X)\to\Db(Y)$ for a proper scheme morphism $f\colon X \to Y$.

We depict divisors consisting of rational curves by their dual intersection graphs: \intext{\tikzchain{1}} denotes a reduced $(-1)$-curve, \intext{\tikzchainmult{3}{2}} a double $(-3)$-curve, \intext{\tikzchain{2,2}} two $(-2)$-curves intersecting in one point, and \intext{\tikzchainnonred{3,3}{2,}} two $(-3)$-curves intersecting in two points.

\section{Properties of divisors: negativity, rigidity, connectivity}
\label{sec:divisor-properties}

\noindent
Let $X$ be a smooth, algebraic, \ie quasi-projective surface, defined over an algebraically closed field $\kk$. A \emph{curve} on $X$ will always mean an irreducible, reduced effective divisor, not necessarily smooth. Therefore, given an effective divisor $D$ on $X$, we will speak of its \emph{curve decomposition} $D=\sum_i c_iC_i$, where the $C_i$ are pairwise distinct curves and all $c_i$ are positive.

In this section, we typographically distinguish whether we \emph{recall} an established notion or \boldemph{introduce} a new one.

\subsubsection*{Intersection numbers and Euler characteristics}

The intersection number of divisors can be computed cohomologically with the Euler characteristic.
Recall that for coherent sheaves $\FF, \mathcal G$ on $X$ with proper support, their Euler characteristic is
$\chi(\FF,\mathcal G) = \sum_i (-1)^i \dim\Ext^i(\FF,\mathcal G)$.

\begin{lemma} 
\label{lem:intersection-numbers}
Let $A$ and $B$ be two effective divisors on $X$. Then
\[ A.B = \chi(\OO_A) + \chi(\OO_B) - \chi(\OO_{A+B}) = - \chi(\OO_A,\OO_B) 
, \]
and, in particular, $A^2 = -\chi(\OO_A,\OO_A)$.
\end{lemma}

\begin{proof}
We start with the well-known formula, see \eg \cite[Ex.~V.1.1]{Hartshorne}: $A.B = \chi(\OO) - \chi(\OO(-A)) - \chi(\OO(-B)) + \chi(\OO(-A-B))$; replacing $X$ by an auxiliary compactification, if needed. This statement holds for all (projective) divisors, not only effective. The first equation follows immediately, using ideal sheaf sequences.
For the second one, compute 
\begin{align*}
  \hspace{-1em}
  \chi(\OO_A,\OO_B) &= \chi(\OO,\OO_B) - \chi(\OO(-A),\OO_B)
                     = \chi(\OO_B) - \chi(\OO_B(A)) = \\
                   &= \big(\chi(\OO)-\chi(\OO(-B))\big) - \big(\chi(\OO(A))-\chi(\OO(A-B))\big)
                      = -A.B \qedhere
\end{align*}
\end{proof}

For an effective divisor $D$ on $X$, we have the Riemann--Roch formula for the Euler characteristic of its structure sheaf:
\[ \chi(\OO_D) = -\tfrac{1}{2} (D^2+D.K_X) . \]
%
%

\begin{lemma}[Decomposition sequence]
If $A+B$ is a sum of two effective divisors, then there is a short exact sequence of sheaves on $X$
\[
0 \to \OO_A(-B) \to \OO_{A+B} \to \OO_B \to 0 .
\]
\end{lemma}

\begin{proof}
We just include the proof for the homological fun of it. Consider the following commutative diagram of ideal sheaf sequences
\[
\begin{tikzcd}[column sep=1.2em, row sep=2.4ex]
   0 \ar[r] & \OO(-A-B) \ar[r] \ar[d, "\iota"] & \OO \ar[r] \ar[d, "\id"] & \OO_{A+B} \ar[r] \ar[d, "\pi"] & 0\\
   0 \ar[r] & \OO(-B)   \ar[r] & \OO \ar[r] & \OO_{B} \ar[r] & 0\\
\end{tikzcd}
\]
where $\iota$ is the inclusion and $\pi$ the restriction.
The Snake Lemma tells us that $\pi$ is surjective with kernel isomorphic to $\coker \iota = \OO_A(-B)$.
\end{proof}

\subsection{Negativity properties of divisors} \label{sec:negativity-properties}

The negativity of a divisor $D$ can be measured in several ways:
\begin{itemize*}
\item $D$ is \emph{negative} if $D^2<0$.
\item $D$ is \emph{negative definite} if $A^2<0$ for all $0 \neq A$ with $\supp(A)\subseteq\supp(D)$.
\item $D$ is \boldemph{negatively closed} if $A^2<0$ for all $0 \led A\leqd D$.
\item $D$ is \boldemph{negatively filtered} if $D = C_1 + \cdots + C_m$ with curves $C_i$ such that $C_i.(C_i+\cdots+C_m) < 0$ for all $i=1,\ldots,m$.
\end{itemize*}
In words, $D$ is negatively closed if all effective subdivisors are negative. By contrast, $D$ is negative definite if every divisor supported on (a subset of) $D$ is negative, \ie the intersection matrix of the curve components of $D$ is negative definite; in particular, $D$ is then negatively closed. Note that negative definiteness depends only on the curve configuration underlying $D$. 

\begin{remark} \label{rem:negative-filtration}
To explain negative filtrations, let $C\led D$ be a curve with $C.D<0$, and set $D' \coloneqq D-C$. The cohomology of the decomposition sequence $\OO_{D'}(D') \into \OO_D(D) \onto \OO_C(D)$ gives $H^0(\OO_{D'}(D')) \cong H^0(\OO_D(D))$, since $H^0(\OO_{C}(D)) = 0$ by the assumption $C.D<0$.

Therefore, by inductively applying this argument to a negative filtration $D = C_1+\cdots+C_m$, we find $H^0(\OO_D(D)) = \cdots = H^0(\OO_{C_m}(C_m)) = 0$, where the last equation uses $C_m^2<0$ from the final step of the filtration. Hence $D$ is rigid as a subscheme in the parlance of \autoref{sec:rigidity-properties}.
\end{remark}

\begin{lemma} \label{lem:negatively-closed-implies-negatively-filtered}
A negatively closed divisor $D$ is negatively filtered.
\end{lemma}

\begin{proof}
More generally, we show that any curve $C$ in a negatively closed divisor $D$ with $C.D<0$ can be extended to a negative filtration of $D$.

First, we note that a negatively closed divisor $D$ contains a curve $C_1 \coloneqq C$ with $C.D<0$. Otherwise, writing $D = \sum n_i C_i$ as a sum of curves yields a contradiction $D^2 = \sum n_i C_i.D \geq 0$. Now proceed by induction with $D-C$ which is negatively closed as well.
\end{proof}

\begin{example} \label{ex:negatively-filtered_not_negatively-closed}
Let $D = 2E + A+A'+B = \tikztreemult{3,3,2}{1}{1,1,1}{2}$. Then $A,A',E,B,E$ is a negative filtration for $D$. However, $D^2 = -4-3-3-2 + 2(2+2+2) = 0$, so that $D$ is not negatively closed. 
\end{example}

\noindent
\autoref{ex:spherelike-not-pullback} gives a negatively closed divisor which is not negative definite.

\subsection{Rigidity properties of divisors} \label{sec:rigidity-properties}

An effective divisor $D$ is called
\begin{itemize*}
\item \emph{rigid} (as a sheaf) if the first order infinitesimal deformations of its structure sheaf $\OO_D$, as a sheaf on $X$, are trivial, \ie $\Ext^1_X(\OO_D,\OO_D)=0$.
\item \emph{rigid as a subscheme} if all first order infinitesimal deformations of $D$ as a closed subscheme of $X$ are trivial. As these infinitesimal deformations are classified by global sections of the normal bundle of $D$ in $X$, this condition amounts to $H^0(\OO_D(D))=0$.
\item \emph{Jacobi rigid} if there are no first order infinitesimal deformations of  $\OO_D$ as a degree zero line bundle on $D$. This amounts to $H^1(\OO_D)=0$ which is the tangent space of the Jacobian $\Pic^0(D)$ at $[\OO_D]$.
\item \emph{rational} if all curves in $D$ are smooth, rational curves.  
\end{itemize*}

\begin{lemma} \label{lem:deformations}
  For any effective divisor $D$ on $X$, there are canonical isomorphisms of vector spaces
\begin{align*}
  \Hom  (\OO_D,\OO_D)  &=  H^0(\OO_D), \\
  \Ext^1(\OO_D,\OO_D)  &=  H^0(\OO_D(D)) \oplus H^1(\OO_D) , \\
  \Ext^2(\OO_D,\OO_D)  &=  H^1(\OO_D(D)) .
\end{align*}
In particular, $D$ rigid $\iff$ $D$ rigid as a subscheme and $D$ Jacobi rigid.
\end{lemma}

\begin{proof}
Denote by $i\colon D\embed X$ the inclusion. We make use of the adjunction $i^* \dashv i_*$ of exact functors between the derived categories $\Db(D)$ and $\Db(X)$. 
Note that $i_*$ coincides with the underived push-forward of sheaves, in contrast to $i^*$. 
Note that there is an isomorphism $\OO_D = i_*\OO_D \cong [\OO(-D) \to \OO]$ in $\Db(X)$ and therefore a decomposition $i^*i_*\OO_D = \OO_D\oplus\OO_D(-D)[1]$. Now we compute
\begin{align*}
  \Hom^\bullet_X(\OO_D,\OO_D) &= \Hom^\bullet_X(i_*\OO_D,i_*\OO_D) \\
                       &= \Hom^\bullet_D(i^*i_*\OO_D,\OO_D) \\
                       &= \Hom^\bullet_D(\OO_D\oplus\OO_D(-D)[1],\OO_D) \\
                       &= \Hom^\bullet_D(\OO_D,\OO_D) \oplus \Hom^\bullet_D(\OO_D(-D)[1],\OO_D) \\
                       &= H^\bullet(\OO_D) \oplus H^\bullet(\OO_D(D))[-1].
\qedhere
\end{align*}
\end{proof}

\begin{example} \label{ex:negative-elliptic-curve}
We give a standard example of a curve that is rigid as a subscheme but not Jacobi rigid: Let $C\subset\IP^2$ be a smooth cubic; in particular, $C^2=9$. Blowing up $\IP^2$ in eleven points lying on $C$, the total transform of $C$ is $D=\bar C+E_1+\cdots+E_{11}$, where $\bar C$ is the strict transform of $C$. As the intersection pairing is preserved, we get $9=C^2=D^2=\bar C^2-11+20$, hence $\bar C$ is an elliptic curve with $\bar C^2=-2$. It is rigid as a subscheme but $\OO_{\bar C}$ deforms in view of $H^1(\OO_{\bar C})=\kk$. By Riemann--Roch, $\bar C.K_X=2$.
\end{example}

\subsection{Connectivity properties of divisors} \label{sec:connectivity-properties}

Let $D$ be an effective divisor.
\begin{itemize*}
\item $D$ is \emph{connected} if it has connected support.
\item $D$ is \emph{well-connected} if $H^0(\OO_D) = \Hom(\OO_D,\OO_D)=\kk$.
\item $D$ is \emph{(numerically) 1-connected} if $A.(D-A) \geq 1$ for any $0 \led A \led D$.
\end{itemize*}
A well-connected divisor is obviously connected. For reduced divisors, the three notions are equivalent.
\begin{example} \label{ex:negative_not_well-connected}
Let $E$ be a $(-1)$-curve, \ie a smooth, rational curve with $E^2=-1$. The decomposition sequence for $2E$ easily gives $H^0(\OO_{2E})=\kk^3$.
\end{example}

The following result of Ramanujan is classical.

\begin{lemma}[{e.g. \cite[Lem. 3.11]{Reid}}] \label{lem:1-connected-is-well-connected}
1-connected divisors are well-connected.
\end{lemma}

\begin{example} \label{ex:well-connected_not_1-connected}
Consider the divisor $D=
\raisebox{-0.5ex}{
  $\tikz
  {
  \node[curve,double] (1) at (0,0)     {\printnum{1}};
  \node[curve,double] (2) at (1.2,0)   {\printnum{2}};
  \node[curve,double] (3) at (0.6,0.4) {\printnum{3}};
  \draw (1) -- (2) -- (3) -- (1);
  }$
}$
consisting of a triangle of rational curves, each of multiplicity two.
The divisor is not 1-connected: $D\red^2=-1-2-3+2\cdot 3=0$. Starting from the decomposition sequence for $D = D\red + D\red$, one can check that $D$ is well-connected.
%
\end{example}

\begin{question}
Is a well-connected divisor whose dual intersection graph is a tree automatically 1-connected?

This is about the reverse of 1-connected $\Longrightarrow$ well-connected.
The previous \autoref{ex:well-connected_not_1-connected} is a counterexample, but it is not a tree.
\end{question}

\noindent
We consider two more conditions related to connectivity:
\begin{itemize*}
\item $D$ is a \emph{tree} if the dual intersection graph of $D$ is a tree and all intersections of irreducible curves in $D$ are transversal.
\item $D$ is \boldemph{1-decomposable} if it can be written as $D=C_1+\cdots+C_m$ where the $C_i$ are irreducible curves which satisfy $C_i.(C_{i+1}+\cdots+C_m)=1$ for all $i=1,\ldots,m-1$.
\end{itemize*}

\begin{remark} \label{rem:1-decompositions}
To explain the definition of 1-decompositions, let $C\led D$ be a \emph{rational} curve with $C.D'=1$, where $D' \coloneqq D-C$. The cohomology of the decomposition sequence $\OO_C(-D') \into \OO_D \onto \OO_{D'}$ gives $H^\bullet(\OO_D) \cong H^\bullet(\OO_{D'})$, since $H^\bullet(\OO_C(-D')) = 0$ by the assumption $C.D'=1$.

Therefore, by applying this argument inductively, a  1-decomposition of a rational divisor $D=C_1+\cdots+C_m$ yields $H^\bullet(\OO_D) = \cdots = H^\bullet(\OO_{C_m}) = \kk$.
\end{remark}

Regarding the definition of trees, many sources include the condition that all curves are smooth and rational. For our purposes it is better to distinguish between these properties (\eg in the next paragraph).

Any divisor that is a reduced tree of curves has a $1$-decomposition, by inductively pruning the leaves. 
In \autoref{lem:1-decomposable-tree} we show that 1-decomposable divisors are trees, and we conjecture that they are 1-connected as well.

\subsection{\boldmath{$(-n)$}-divisors}

In this article, we are interested in well-connected, rigid divisors. Concretely, for an effective divisor $D$ this means $H^0(\OO_D) = \kk$ and $\Ext^1(\OO_D,\OO_D) = 0$. 
Recall $\Ext^1(\OO_D,\OO_D) = H^0(\OO_D(D)) \oplus H^1(\OO_D)$.

\begin{remark}
$\OO_D\in\Db(X)$ is thus an object of a triangulated category with $\Hom^0(\OO_D,\OO_D)=\kk$ and $\Hom^1(\OO_D,\OO_D)=0$. In topological parlance, \eg the category of spectra, one might call $\OO_D$ a `simply connected' object. We prefer to stick with `well-connected rigid' to stress the rigidity of $D$.
\end{remark}

Often, it is convenient to specify $\ext^2(\OO_D,\OO_D)$.
For well-connected, rigid $D$, we have
 $h^1(\OO_D(D)) = \ext^2(\OO_D,\OO_D) = \chi(\OO_D,\OO_D)-1 = -D^2-1$.
Thus:
\begin{itemize*}
\item $D$ is a \boldemph{\boldmath$(-n)$-divisor} if $D$ is well-connected and rigid with $D^2=-n <0$.
\end{itemize*}

By \autoref{lem:deformations}, $D$ is a $(-n)$-divisor if and only if $H^1(\OO_D(D))=\kk^{n-1}$ and $H^0(\OO_D)=\kk$ and $H^1(\OO_D)=H^0(\OO_D(D))=0$. Assuming the other three provisions, $H^1(\OO_D(D))=\kk^{n-1}$ can be replaced with $D^2=-n$.

We mention three important special cases of this notion:
\begin{itemize*}
\item $D$ is \emph{exceptional} if $D$ is a $(-1)$-divisor.
\item $D$ is \boldemph{spherelike} if $D$ is a $(-2)$-divisor.
\item $D$ is \emph{spherical} if it is spherelike and invariant under twisting with the canonical bundle of $X$, \ie $\OO(K_X)|_D \cong \OO_D$.
\end{itemize*}
Saying that $D$ is exceptional, or spherelike, or spherical, is just a shortcut for the sheaf $\OO_D \in \Db(X)$ being an exceptional object (\cite[\S 8.3]{Huybrechts}), or a 2-spherelike object (\cite{HKP}), or a 2-spherical object (\cite[\S 8.1]{Huybrechts}), respectively.

In the definition of $(-n)$-divisors, the condition of being well-connected is crucial. The next remark spells out what happens when it is dropped, in the first non-trivial case ($n=2$).

\begin{remark} Let $D$ be a rigid divisor with $D^2=-2$, \ie we drop the condition that $D$ be well-connected from the definition of spherelike. Due to
  $1 \leq \hom(\OO_D,\OO_D) \leq \chi(\OO_D,\OO_D) = -D^2 = 2$,
there are only two options:
\begin{itemize}
\item $\hom(\OO_D,\OO_D)=1$, \ie $D$ is spherelike as defined above;
\item $\hom(\OO_D,\OO_D)=2$, 
      \ie $\OO_D \in \Db(X)$ a 0-spherelike object; see \cite{HKP}.
\end{itemize}
Because $\kk$ is algebraically closed, the 0-spherelike case in turn leads to either $\End(\OO_D) = \kk[x]/x^2$ or $\End(\OO_D)=\kk\times\kk$. Examples for these two types are \intext{\tikzchainmult{1,2}{2,1}} and the disconnected divisor \intext{\tikzchainnonred{1,1}{0,}}, respectively. The first of these shows that a merely connected and rigid divisor is not necessarily a $(-n)$-divisor, one has to ask for well-connectedness.
%
\end{remark}

\section{Relations among divisors properties}

\subsection{Properties of rigid divisors}

\begin{lemma} \label{lem:curves-rigid=negative}
A rigid divisor is negative.
\end{lemma}

\begin{proof}
We have
 $-D^2 = \chi(\OO_{D},\OO_{D}) = \hom(\OO_D,\OO_D) + \ext^2(\OO_D,\OO_D) > 0$
with $\Ext^1(\OO_D,\OO_D)=0$ because $\OO_{D}$ is rigid.
\end{proof}



\begin{lemma} \label{lem:rigid-descends}
Any effective subdivisor of a rigid divisor is rigid.

Likewise with `rigid' replaced by `Jacobi rigid' or `rigid as a subscheme'.
\end{lemma}

\begin{proof}
Let $D=A+B$ be a decomposition of $D$ into two effective divisors.

First assume that $D$ is Jacobi rigid, \ie $H^1(\OO_D)=0$. 
The long exact cohomology sequence of the decomposition sequence
 $\OO_B(-A) \into \OO_D \onto \OO_A$
contains a surjection $H^1(\OO_D) \onto H^1(\OO_A)$, thus $H^1(\OO_A)=0$, and $A$ is also Jacobi rigid.

Now assume $D$ is rigid as a subscheme, \ie $H^0(\OO_D(D))=0$.
Tensoring the other decomposition sequence
$\OO_A(-B) \into \OO_D \onto \OO_B$ with $\OO(D)$
yields an injection $H^0(\OO_A(A)) \into H^0(\OO_D(D))$. So $H^0(\OO_A(A))=0$, and hence $A$ is rigid as a subscheme of $X$.

By \autoref{lem:deformations}, both cases together imply that $A$ is rigid when $D$ is.
\end{proof}

These two lemmas combine to the following statement:

\begin{corollary} \label{cor:rigid-totally}
A rigid divisor is negatively closed.
\end{corollary}

\begin{example}
The sum of two rigid divisors does not have to be rigid:
Let $C$ be a curve with $C^2=0$, \eg the fibre of a ruled surface $X$. Then $C$ obviously deforms at least locally. Next, let $\tilde X$ be the blow-up of $X$ in any point of $C$ and denote by $\tilde C = \intext{\tikzchain{1,1}}$ the total transform of $C$. Then $\tilde C^2=0$ and $\tilde C$ is not rigid, as $\Ext^1_X(\OO_{\tilde C},\OO_{\tilde C})=\Ext^1_X(\OO_C,\OO_C)\neq0$.
\end{example}

\begin{lemma} \label{lem:rigid-tree-rat}
A connected, Jacobi rigid divisor is a rational tree.
\end{lemma}

\begin{proof}
Let $D=\sum_i c_iC_i$ be the curve decomposition of $D$, and assume $D$ is connected and Jacobi rigid. By \autoref{lem:rigid-descends}, each curve $C_i$ is Jacobi rigid. Thus $p_a(C_i)=1-\chi(\OO_{C_i})=1-h^0(\OO_{C_i})+h^1(\OO_{C_i})=1-1+0=0$, \ie $C_i$ has arithmetic genus zero and is smooth and rational (see \cite[\S II.11]{BHPV}).

The divisor $D\red=\sum_i C_i$ of $D$ is reduced by definition, connected by assumption and also has arithmetic genus
$1-\chi(\OO_{D\red})=1-1+0=0$ (note that reduced and connected implies $H^0(\OO_{D\red})=\kk$) and, again by \cite[\S II.11]{BHPV}, a tree of smooth, rational curves.
\end{proof}

\begin{example} \label{ex:negative-definite_not_tree}
Let $D = \intext{\tikzchainnonred{2,3}{2,}}$ be a $(-2)$-curve and $(-3)$-curve intersecting transversally in two points. Then $D$ is easily seen to be negative definite. However, $D$ is not a tree and, in particular, not Jacobi rigid.
\end{example}

\begin{example} \label{ex:negative-rational-tree_not_rigid}
Consider $D = 2Z+2A+2B+2C = \tikztreemult{3,3,3}{1}{2,2,2}{2}$ where all curves are rational with $A^2=B^2=C^2=-3$ and the central curve $Z^2=-1$. An easy computation yields that $D$ is negatively closed. Put $D' \coloneqq D_{\red}$. The decomposition sequence for $D = D' + D'$ shows $H^1(\OO_D)=H^1(\OO_{D'}(-D'))$ because $\OO_{D'}$ is rigid and well-connected. The decomposition sequence for $D' = (A+B+C) + Z$, twisted  by $-D'$, is
\[ \OO_{A+B+C}(-Z-D') \into \OO_{D'}(-D') \onto \OO_Z(-D') . \]
Its cohomology sequence ends in
  $H^1(\OO_{D'}(-D')) \onto H^1(\OO_Z(-D')) = \kk$,
because of $Z.(-D') = Z.(-A-B-C-Z) = -2$.
Therefore we have $H^1(\OO_D) = H^1(\OO_{D'}(-D'))\neq0$, and $D$ is a negatively closed divisor which is a rational tree but not Jacobi rigid.
\end{example}

\noindent
See also \autoref{ex:elliptic-sing} for a negative definite divisor that is a tree but not Jacobi rigid.

\subsection{Properties of 1-decomposable divisors}

The results of this subsection are not used for the numerical criterion for $(-n)$-divisors \autoref{thm:numerical-criterion}, and the reader may jump ahead to \autoref{sub:overview}.

\begin{lemma} \label{lem:1-decomposable-tree}
A 1-decomposable divisor is a tree.
\end{lemma}

\begin{proof}
Note that a 1-decomposable divisor is automatically connected, so we only need to show that there are no cycles.
  Let $D=C_1+\cdots+C_m$ be a 1-decomposition of $D$. For a contradiction, assume that $C_{i_1},\ldots,C_{i_k}$ is a cycle, where $i_1<i_2<\ldots<i_k$. Without loss of generality, the curve $C_{i_1}$ does not occur among $C_j$ with $j>i_1$. The inequality
$
  1 =    C_{i_1}.(C_{i_1+1}+C_{i_1+2}+\cdots+C_m)
  \geq C_{i_1}.(C_{i_2}+C_{i_3}+\cdots+C_{i_k}) \geq 2
$
is the desired contradiction.
\end{proof}

\begin{lemma} \label{lem:1-decomposable-negative-curves}
  If $C$ is a multiple curve in a 1-decomposable divisor $D$, then $C^2\leq0$ with the single exception of $D=2C$ and $C^2=1$.
\end{lemma}

\begin{proof}
  For a $1$-decomposition $D=C_1+\cdots+C_m$, let $D' \coloneqq C_{i+1}+\cdots+C_m$ with $i$ minimal such that $C_i=C$. By assumption, $C\leqd D'$.

Denote by $w\geq2$ the multiplicity of $C$ in $D$, so that $C$ has multiplicity $w-1$ in $D'$. Next, denote by $v \coloneqq C.(D'-(w-1)C)$ the valency of $C$ in $D'$, \ie the number of curves in $D'$ (with multiplicities) intersecting $C$. Then $1 = C.D' = (w-1)C^2 + v$ by the definition of 1-decomposition, hence $C^2 = (1-v)/(w-1)$. Thus $C^2\leq 0$ unless $v=0$.

However, if $v=0$ then $C^2=1$ and $w=2$. As $D'$ is 1-decomposable and hence connected, $v=0$ implies $D'=C$. This enforces $D=2C$, lest we get a contradiction $1 = C_{i-1}.(C_i+\cdots+C_m) = C_{i-1}.(C+D') = C_{i-1}.2C$.
\end{proof}

\begin{conjecture} \label{conj:1-decomposable-1-connected}
A 1-decomposable divisor is 1-connected.
\end{conjecture}

\begin{remark}
This conjecture could be more generally seen as a purely combinatorial statement about finite weighted graphs. 
Note that there are graphs which cannot be realised on any surface as a divisor, \eg \autoref{ex:violates-HIT} is a reduced tree (hence 1-decomposable).
For dual graphs of rational divisors, it follows from \autoref{lem:1-decomposition-given}, \autoref{lem:curvelike-num1conn} and \autoref{lem:has-1-decomposition}.
Hence the conjecture can be reduced to show:
If $D$ is a 1-decomposable graph on $X$, then there exists a 1-decomposable, rational graph $D'$ on some surface $X'$.
\end{remark}

\subsection{About divisor properties} \label{sub:overview}

The multitude of divisor properties can be overwhelming. In this short section, we attempt to categorise them.

\begin{definition} \label{def:divisor-properties}
Let $(P)$ be a property of effective divisors. $(P)$ is called
  
\begin{itemize}
\item \emph{closed} if $0 \led D' \led D$ and $D$ satisfies $(P)$,
      then $D'$ satisfies $(P)$;
\item \emph{open} if $D \led D''$ with $\supp(D)=\supp(D'')$ and $D$ satisfies $(P)$,
      then $D''$ satisfies $(P)$;
\item \emph{numerical} if $D = c_1C_1+\cdots+c_mC_m$ is the curve
      decomposition of a divisor satisfying $(P)$ and $D'=c_1C_1'+\cdots+c_mC_m'$
      with $C_i.C_j=C'_i.C'_j$ for all $i,j$, then $D'$ satisfies $(P)$;
\item \emph{birational} if it is preserved under contractions and blow-ups, \ie for a blow-up $\pi\colon X'\to X$ in a point, $D \subset X$ satisfies $(P)$ if and only if $\pi^*D \subset X'$ satisfies $(P)$;
\item \emph{homological} if it depends only on the graded vector space $\Ext^\bullet_X(\OO_D,\OO_D)$.
\end{itemize}

\end{definition}

\begin{remark}
The topological terminology comes from viewing the set of effective divisors on $X$ as a poset under $\leqd$, and giving it the Alexandrov topology. The notion `open' only makes sense if, given a divisor $D$, we restrict to the subspace of divisors supported on $D$.
\end{remark}

\begin{example} 
Rigid, Jacobi rigid, rigid as subscheme (by \autoref{lem:rigid-descends}) and negatively closed (by definition) are closed properties. 

A property is open and closed if it depends only on the reduced divisor, \ie the underlying curve configuration. Examples are negative definite (by definition) and rigid on exceptional loci of rational singularities.

Homological properties are well-connected, rigid, $(-n)$-divisor and spherelike. A homological property is birational.
\end{example}

\begin{remark}
A property $(P)$ is numerical if it depends on the intersection matrix and the coefficients of the curves making up $D$ (up to permutation). This is different from asking that $(P)$ only depends on the numerical class $[D]\in \mathrm{NS}(X)$ as the matrix and the coefficients do not determine the classes $[C]$ of the curves in $D$.
Equivalently, $(P)$ can be checked on the weighted (by multiplicities) dual intersection graph of $D$.
By \autoref{thm:numerical-criterion}, the property being a $(-n)$-divisor is numerical when restricted to rational divisors.
\end{remark}

For a naturally occurring property which is numerical, but not birational, see \autoref{rem:nonbirational-property}. By contrast, the next statement shows that all properties introduced in \autoref{sec:divisor-properties} are birational.

\begin{proposition} \label{prop:birational-properties}
The following properties of effective divisors on surfaces are birational: negative definite, negatively closed, negatively filtered, $(-n)$-divisor, rigid, rigid as a subscheme, Jacobi rigid, 1-decomposable, 1-connected, well-connected, connected, rational.
\end{proposition}

\begin{proof}
We can assume that $\pi\colon \tilde X \to X$ is the blow-up in a single point, with an effective divisor $D$ on $X$ and its total transform $\tilde D = \pi^* D$ on $\tilde X$. It is sufficient to show that a property holds for $D$ if and only if it holds for $\pi^* D$.
  
The claim is obvious for these two properties: connected and rational.

For the homological notions (rigid, well-connected, $(-n)$-divisor), use that $\pi^*$ is fully faithful:
  $\Ext^i(\OO_D,\OO_D) = \Ext^i(\pi^*\OO_D,\pi^*\OO_D) = \Ext^i(\OO_{\tilde D},\OO_{\tilde D})$
for all $i$. Since $\pi^*\OO_X = \OO_{\tilde X}$, we also find that Jacobi rigidity is birational. By \autoref{lem:deformations}, rigidity as a subscheme is thus birational, too.

Now let $D=C_1+\cdots+C_m$ be a 1-decomposition of $D$. Write $C'_i$ for the strict transform of $C_i$. We get a 1-decomposition of $\smash{\tilde D}$ by replacing each $C_i$ with its total transform $\pi^*C_i$, where the exceptional divisor comes first. More explicitly, if the blown-up point is not on $C_i$, then $\smash{\pi^* C_i=C'_i}$ is unchanged; if the point is on $C_i$, then take $\smash{\pi^* C_i=E+C'_i}$. This is again a 1-decomposition by $E^2=-1$ and $\smash{{C'_i}^2=C_i^2-1}$. 
Conversely, a 1-decomposition $\tilde D=C'_1+\cdots+C'_{m'}$ can always be rearranged such that $E$ is in front of a curve $C'_i$ with $C'_i.E=1$ (as the intersection of $E$ with curves untouched by the blow-down is always zero), so we can revert the above process by replacing each $C'_i$ by $C_i$, and dropping all $E$.

A similar procedure applies to negative filtrations.

Regarding negatively closed divisors, we use some standard facts about the intersection products on $X$ and $\tilde X$. For any divisors $A,B$ on $X$, we have $\pi^*A.\pi^*B=A.B$. This immediately shows that $\tilde D$ negatively closed implies $D$ negatively closed. Conversely, for divisors $B$ on $X$ and $\tilde A$ on $\tilde X$ we have
$B.\pi_*\tilde A = \pi^*B.\tilde A$. If $D$ is negatively closed and  $0 \leqd \tilde A \leqd \tilde D$, then we use that fact with $B \coloneqq \pi_*\tilde A \leqd D$, obtaining
\[
  0  >    (\pi_*\tilde A)^2
     =    (\pi^*\pi_*\tilde A).\tilde A
     =    (\tilde A + (\tilde A.E)E).\tilde A
     =    \tilde A^2 + (\tilde A.E)^2
     \geq \tilde A^2 , 
\]
showing that $\tilde D$ is negatively closed.

$D$ is negative definite if all multiples $kD$, with $k>0$ are negatively closed. Hence the last paragraph also proves that negative definite is birational. 

Finally, for 1-connectedness, similar reasoning works as for negatively closed. If $\tilde D$ is 1-connected, then $D$ is, too: $A.(D-A)=\pi^*A.\pi^*(D-A)=\tilde A.(\tilde D-\tilde A)\geq 1$. Conversely, if $D$ is 1-connected and $0 \leqd \tilde A\leqd \tilde D$, then
\begin{align*}
  1 &\leq \pi_*\tilde A.(D-\pi_*\tilde A)
     =    \pi_*\tilde A.\pi_*(\tilde D-\tilde A)
     =    \pi^*\pi_*\tilde A.(\tilde D-\tilde A) \\
    &=    (\tilde A + (\tilde A.E)E).(\tilde D-\tilde A)
     =    \tilde A.(\tilde D-\tilde A) - (\tilde A.E)^2 ,
\end{align*}
using $\tilde D.E=\pi^*D.E=0$. Hence $\tilde A.(\tilde D-\tilde A) \geq 1+(\tilde A.E)^2 \geq 1$.
\end{proof}

\section{Numerical characterisation of $(-n)$-divisors}

\noindent
In this section, we prove a characterisation of well-connected rigid divisors in terms of intersection numbers. This is useful because it turns a homologically defined notion into something much easier to test in practice. The criterion is about the existence of particularly nice curve decompositions; in examples, this is the easiest way to check that a divisor is a $(-n)$-divisor.

For the convenience of the reader, we recall the notions appearing in the theorem: let $D = C_1+\cdots+C_m$ be a curve decomposition of an effective divisor. The sequence of curves $C_1,\ldots,C_m$ is called a
\begin{itemize}
\item \emph{$1$-decomposition} if  $C_i.(C_{i+1}+\cdots+C_m)=1$ for all $i=1,\ldots,m-1$;
\item \emph{negative filtration} if $C_i.(C_i+\cdots+C_m)<0$ for all $i=1,\ldots,m$.
\end{itemize}

\begin{theorem} \label{thm:numerical-criterion}
For an effective divisor $D$, these conditions are equivalent:
\begin{enumerate}
\item $D$ is well-connected and rigid, \ie a $(-n)$-divisor for some $n>0$.
\item $D$ is rational, $1$-decomposable and negatively filtered.
\end{enumerate}
\end{theorem}

\begin{corollary} \label{cor:numerical-criterion}
An effective rational divisor $D$ is a $(-n)$-divisor if and only if $D$ is $1$-decomposable and negatively filtered with $D.K_X=n-2$.  
\end{corollary}

\begin{proof}
This follows from Riemann--Roch: $1 = \chi(\OO_D)=-\frac12(D^2+D.K)$.
\end{proof}

\begin{remark} \label{rem:numerical-criterion}
This corollary works well for computing spherelike examples, where we can check $D.K_X=0$, which is easier to calculate than $D^2$.
  
The conditions `$D$ rational' and `$D.K = n-2$' are not numerical in the sense of \autoref{def:divisor-properties}. This is unavoidable in view of negative non-rational curves, as in \autoref{ex:negative-elliptic-curve}.

There cannot be a purely numerical characterisation of (non-rational) rigid divisors, at least not with the properties of \autoref{sub:overview}. For instance, \autoref{ex:negative-rational-tree_not_rigid} contains a rational, negatively closed tree but non-rigid divisor.
\end{remark}

\begin{example} \label{ex:sph-2-321}
The divisor $D = B + 2C + C' + E =$ \tikztreemult{3,2,1}{2}{1,1,1}{2} such that $B^2=-3, E^2=-1$ and $C^2=C'^2=-2$ is spherelike, \ie a $(-2)$-divisor, by \autoref{cor:numerical-criterion}: $C,B,C',C,E$ is both a 1-decomposition and a negative filtration, and the conditions $D.K=0$ or $D^2=-2$ are easy to check.



\end{example}

The \hyperlink{proof:numerical-criterion}{proof} of \autoref{thm:numerical-criterion} uses a couple of lemmas. 

\begin{lemma} \label{lem:curvelike-num1conn}
A Jacobi rigid, well-connected divisor is $1$-connected.
\end{lemma}

\begin{proof}
Let $0 \led A \led D$ and $B \coloneqq D-A$; we need to show $A.B \geq 1$. Now $D$ is Jacobi rigid, hence $A$ and $B$ are as well, by \autoref{lem:rigid-descends}. Thus $H^1(\OO_A)=H^1(\OO_B)=0$ and specifically $\chi(\OO_A)\geq1$ and $\chi(\OO_B)\geq1$. About $D$ we know $H^0(\OO_D)=\kk$ and $H^1(\OO_D)=0$, so $\chi(\OO_D)=1$. By \autoref{lem:intersection-numbers},
\[
   A.B = \chi(\OO_A) + \chi(\OO_B) - \chi(\OO_D)
      \geq 1 + 1 - 1 = 1 . \qedhere
\]
\end{proof}

\begin{lemma} \label{lem:negatively-filtered-implies-subscheme-rigid}
A negatively filtered divisor is rigid as a subscheme.
\end{lemma}

\begin{proof}
By \autoref{rem:negative-filtration}, if $D = C_1 + \cdots + C_m$ is a negative filtration, then
 $H^0(\OO_D(D)) = \cdots = H^0(\OO_{C_m}(C_m)) = 0$, the last equality from $C_m^2<0$.
\end{proof}

Part of the next statement is a weaker version of \autoref{conj:1-decomposable-1-connected}. 

\begin{lemma} \label{lem:1-decomposition-given}
If $D$ is a 1-decomposable, rational divisor, then $D$ is well-connected and Jacobi rigid, \ie $H^0(\OO_D)=\kk$ and $H^1(\OO_D)=0$.
\end{lemma}

\begin{proof}
This follows from \autoref{rem:1-decompositions}: if $D=C_1+\cdots+C_m$ is a $1$-decomposition, then 
 $H^\bullet(\OO_D) = \cdots = H^\bullet(\OO_{C_m}) = \kk$.
\end{proof}


\begin{lemma} \label{lem:self-squares}
Let $D=C_1+\cdots+C_m$ be a curve decomposition of a rational divisor. Then $-\sum_i C_i^2 = D.K_X + 2m$.
\end{lemma}

\begin{proof}
By Riemann--Roch for $C_i \cong \IP^1$, we have $C_i^2+2 = -C_i.K_X$. Summing these up, we get $\sum C_i^2 + 2m = -D.K_X$.
\end{proof}

\begin{lemma} \label{lem:has-1-decomposition}
A Jacobi rigid, $1$-connected divisor is $1$-decomposable.
\end{lemma}

\begin{proof}
By \autoref{lem:rigid-tree-rat}, a Jacobi rigid divisor is rational.

Given a $1$-connected divisor $D$, we have $C.(D-C)\geq1$ for all curves $C$ in $D$. We first claim that there actually is some $C$ with equality. So assuming that $C.(D-C)\geq 2$ throughout, we get an inequality by summing over all curves $C \led D$ and invoking \autoref{lem:self-squares} and Riemann--Roch:
\[
  2m  \leq \sum_{\mathclap{C \led D}} C.(D-C)
       =   \sum_{\mathclap{C \led D}} C.D - \sum_{\mathclap{C \led D}} C^2
       =   D^2 + D.K + 2m
       =   -2\chi(\OO_D) + 2m .
\]
This achieves the desired contradiction, because $\chi(\OO_D)=h^0(\OO_D)>0$ from $D$ Jacobi rigid. So there is a curve $C_1$ such that $C_1.(D-C_1)=1$.

For the induction step, observe that Jacobi rigidity is passed from $D$ on to $D' = D-C_1$ by \autoref{lem:rigid-descends}. The only ingredient missing for the above argument is that $D'$ stays $1$-connected. 
To see this, assume $D = D' + C_1$ is 1-connected with $C_1.D'=1$ and let $D' = A+(D'-A)$ with $0 \led A \led D'$. Then we can write $D'+C_1$ as a sum of effective subdivisors in two ways: $D'+C_1 = (A+C_1)+(D'-A) = A+(D'-A+C_1)$. Thus
\[ \begin{array}{rclcl}
  1 &\leq& (A+C_1).(D'-A) &=& A.(D'-A) -C_1.A  + C_1.D' , \\[0.75ex]
  1 &\leq& A.(D'-A+C_1)   &=& A.(D'-A) +C_1.A.
\end{array} \]
We deduce $\frac{1}{2}\leq A.(D'-A)$, so that $D'$ is 1-connected, as claimed.
\end{proof}

\begin{remark}
\label{rem:can-complete-1decomps}
The proof of this lemma also shows that any partial $1$-decomposition $D=C_1+\cdots+C_l+D'$ of a Jacobi rigid, $1$-connected divisor can be completed to a $1$-decomposition.
\end{remark}

\begin{proof}[Proof of \autoref{thm:numerical-criterion}] \hypertarget{proof:numerical-criterion}{}
$(1)\implies(2)$: $D$ is a $(-n)$-divisor for some $n>0$, so by definition well-connected and rigid. Thus $D$ is $1$-connected by \autoref{lem:curvelike-num1conn}, and rational by \autoref{lem:rigid-tree-rat}. Then $D$ has a $1$-decomposition by \autoref{lem:has-1-decomposition}.
As $D$ is rigid, it is negatively filtered by \autoref{cor:rigid-totally} and \autoref{lem:negatively-closed-implies-negatively-filtered}.

$(2)\implies(1)$: $D$ is rigid and well-connected by \autoref{lem:negatively-filtered-implies-subscheme-rigid} and \autoref{lem:1-decomposition-given}.
\end{proof}

%
%
%
%
%

\section{Minimal $(-n)$-divisors}
\label{sec:minimally-curvelike}

\subsection{Modifying \boldmath{$(-n)$}-divisors: blow-up and blow-down}
\label{sub:modify-blowup}

Consider a blow-up $\pi \colon \tilde X \to X$ in a point $P\in X$, with exceptional $(-1)$-curve $E$, and let $D$ be a $(-n)$-divisor on $X$ for some $n>0$. As $\pi^* \colon \Db(X) \to \Db(\tilde X)$ is fully faithful, $\Hom^\bullet(\pi^*\OO_D, \pi^*\OO_D) = \Hom^\bullet(\OO_D,\OO_D)$. Hence $\pi^* D$ stays a $(-n)$-divisor.

If $P \notin \supp(D)$, then $\pi^* D$ is the same curve configuration as $D$.
Otherwise, denote by $C_1,\ldots,C_k$ the curves in $D$ containing $P$, and let $n_1,\ldots,n_k$ be their multiplicities in $D$, respectively. Then the exceptional curve $E$ appears in $\pi^* D$ with multiplicity $n_1+\cdots+n_k$, and the self-intersection numbers of the strict transforms are $\tilde C_i^2 = C_i^2-1$, as $C_i$ is smooth and contains $P$. For later reference, we record the contraction criterion in the next proposition. In that situation, we say that the $(-n)$-divisor can be blown down.

\begin{proposition} \label{prop:contracting-curves}
Let $D$ be a $(-n)$-divisor on $X$ and $E\led D$ a $(-1)$-curve with $E.D=0$. Let $\pi\colon X \to X'$ be the contraction of $E$. Then $\pi(D)$ is a $(-n)$-divisor on $X'$.
\end{proposition}

\subsection{Modifying \boldmath{$(-n)$}-divisors: twisting spherical divisors}
\label{sub:modify-twisting}

In a certain analogy to the contraction of $(-1)$-curves from a $(-n)$-divisor, there is a construction for $(-2)$-curves, using spherical twists. We give the definitions first, and then proceed to explain why they make sense.

\begin{definition} \label{def:spherelike-component}
 A \emph{spherelike component} of a $(-n)$-divisor $D$ is a spherelike subdivisor $A \led D$ such that $D-A$ is a $(-n)$-divisor.

$A$ is a \emph{spherical component} of $D$ if it is a spherical divisor and a spherelike component of $D$. In this case, $D$ is said to be obtained by \emph{twisting $A$ on to $D-A$}; likewise, $D-A$ is said to be obtained by \emph{twisting $A$ off $D$}. 
\end{definition}

There is a simple check when a $(-2)$-curve is a spherelike, and hence spherical, component:

\begin{lemma} \label{lem:twistonto-single-curve}
If $D$ is a $(-n)$-divisor and $C$ a $(-m)$-curve with $m\geq 2$ and $D.C=1$, then $D+C$ is a $(-n-m+2)$-divisor.

In particular, if $C$ is a $(-2)$-curve, then $D+C$ stays a $(-n)$-divisor.
\end{lemma}

\begin{proof}
This follows from the decomposition sequences
\[
\OO_C(-D) \into \OO_{D+C} \onto \OO_D   \quad\text{and}\quad
\OO_D(D) \into \OO_{D+C}(D+C) \onto \OO_C(D+C)
\]
by taking cohomology and using $C.(-D) = -1$ and $C.(D+C) = 1-m \leq -1$, respectively.
\end{proof}

We turn to some easy and general properties of spherelike components:

\begin{lemma} \label{lem:spherelike-decomp}
Let $A \led D$ be a spherelike component. Then
\begin{enumerate}
\item $A.(D-A) = 1$ and $H^\bullet(\OO_A(A-D))=0$,
\item $\Hom^\bullet(\OO_{D-A},\OO_A(A-D)) = \kk[-1]$,
\end{enumerate}
\end{lemma}

\begin{proof}
(1) Write $B \coloneqq D-A$. Then
  $-n = D^2 = (A + B)^2 = A^2 + 2A.B + B^2 = -2 + 2A.B -n$,
since $A$ is a spherelike component, and hence $A.B=1$.
Next, the decomposition sequence yields
  $H^\bullet(\OO_A(-B)) \to H^\bullet(\OO_D) \to H^\bullet(\OO_B)$,
and by $H^\bullet(\OO_D) \cong H^\bullet(\OO_{B}) \cong \kk$, the second map is an isomorphism induced by restriction of global sections. Hence $H^\bullet(\OO_A(A-D))=0$.

(2) We turn to $\Hom^\bullet(\OO_B,\OO_A(-B))$. For this, apply $\Hom^\bullet(\blank,\OO_A(-B))$ to $\OO(-B) \into \OO \onto \OO_B$ and get
\[
H^\bullet(\OO_A) \from H^\bullet(\OO_A(-B)) \from \Hom^\bullet(\OO_B,\OO_A(-B)).
\]
The left-hand term is isomorphic to $\kk$, and the middle term vanishes by (1), so the right-hand term is isomorphic to $\kk[-1]$, \ie $\Ext^1(\OO_B,\OO_A(-B)) = \kk$ and $\Hom(\OO_B,\OO_A(-B)) = \Ext^2(\OO_B,\OO_A(-B)) = 0$.
\end{proof}

We justify the terminology of \autoref{def:spherelike-component}: for any $F \in \Db(X)$, the \emph{twist functor} $\TTT_F\colon \Db(X) \to \Db(X)$ is defined on objects $A$ by the triangles
 $\Hom^\bullet(F,A) \otimes F \to A \to \TTT_F(A)$,
\ie $\TTT_F(A)$ is the cone of the canonical evaluation morphism. See \cite[\S 2.1]{HKP} for how these cones become functorial in an appropriate setting, \eg for $\Db(X)$.
By \autoref{prop:spherelike-decomposition-with-one-spherical}, spherical components can be twisted off $(-n)$-divisors, leaving smaller $(-n)$-divisors.

\begin{proposition}[{\cite[\S 8.1]{Huybrechts}, \cite[Lem.~3.1]{HKP}}]
Given two divisors $D, D'$ with $D$ effective, the functor $\TTT_{\OO_D(D')}$ is an autoequivalence of $\Db(X)$ if and only if $D$ is spherical.
In this case, $\TTT_{\OO_D(D')}$ is called \emph{spherical twist functor}.
\end{proposition}

\begin{proposition} \label{prop:spherelike-decomposition-with-one-spherical}
Let $A \led D$ be a spherical component. Then
\begin{enumerate}
\item $\Hom^\bullet(\OO_A(A-D),\OO_D) = \kk$, and
\item $\TTT_{\OO_A(A-D)}(\OO_D) \cong \OO_{D-A}$. 
\end{enumerate}
\end{proposition}

\begin{proof}
We apply $\Hom^\bullet(\blank,\OO_A)$ to $\OO(-A) \into \OO(B) \onto \OO_D(B)$ where $B=D-A$, obtaining the triangle 
\[
H^\bullet(\OO_A(A)) \from H^\bullet(\OO_A(-B)) \from\Hom^\bullet(\OO_D,\OO_A(-B)) .
\]
The middle term vanishes by \autoref{lem:spherelike-decomp}(1). Now $A$ is a spherical divisor, \ie $H^\bullet(\OO_A(A)) = \kk[-1]$ and $\OO_A(K) \cong \OO_A$, and Serre duality yields (1):
\[ 
\kk \cong H^\bullet(\OO_A(A))[1] \cong \Hom^\bullet(\OO_D,\OO_A(-B))[2] \cong \Hom^\bullet(\OO_A(-B),\OO_D)^* .
\]
Now (2) follows, since $\TTT_{\OO_A(-B)}(\OO_D)$ is defined by the triangle
\[
\Hom^\bullet (\OO_A(-B),\OO_D) \otimes \OO_A(-B) \to \OO_D \to \TTT_{\OO_A(-B)}(\OO_D) ,
\]
which reduces to the short exact sequence $\OO_A(-B) \into \OO_D \onto \OO_B$.
\end{proof}

\begin{corollary}
\label{cor:curvelike-2-spherical}
Let $D$ be a well-connected rigid divisor.
Then $D$ is spherical precisely if $D$ consists entirely of $(-2)$-curves.
\end{corollary}

\begin{proof}
$(\Longrightarrow)$
Assume that $D$ is spherical and let $C\leqd D$ be any curve. The restriction $\OO_D\onto\OO_C$ induces $\OO_D = \OO_D\otimes\omega_X \onto \OO_C\otimes\omega_X = \OO_C(C.K_X)$. Hence $C.K_X\geq0$ and, summing up,
$0 \leq \sum_{C\leqd D}C.K_X = D.K_X = 0$ by \autoref{lem:self-squares}. 
Therefore $C.K_X=0$ and each $C\leqd D$ is a $(-2)$-curve.

$(\Longleftarrow)$ Let $D$ be well-connected and rigid, consisting of $(-2)$-curves. The claim holds if $D$ is a single curve. We proceed inductively: by \autoref{thm:numerical-criterion}, $D$ is 1-decomposable; let $C\leqd D$ be a curve with $C.D=1$ and assume that $D-C$ is spherical.
By \autoref{lem:twistonto-single-curve}, $D$ is spherelike.
%
The decomposition sequence for $D$, tensored with $\OO(K_X)$, is
 $\OO_C(C-D) \into \OO_{D}(K_X) \onto \OO_{D-C}$,
as both $\OO_{D-C}$ and $\OO_C(C-D)$ are invariant by induction.
By \autoref{lem:spherelike-decomp}(2), $\Ext^1(\OO_{D-C},\OO_C(C-D)) \cong \kk$, and thus $\OO_{D}(K_X) \cong \OO_{D}$.
\end{proof}

\begin{remark}
Consider an abstract setting: if $S_1,\ldots,S_n$ are $d$-Calabi--Yau ($d$-CY) objects in a triangulated category $\DD$, \ie $\Hom(S_i,\blank)\cong\Hom(\blank,S_i[d])^*$, will an object generated by $S_1,\ldots,S_n$ be $d$-CY as well?

\autoref{cor:curvelike-2-spherical} provides a very partial positive answer for a tree of $(-2)$-curves $D = C_1+\cdots+C_n$. Then $\OO_D \in \Db(X)$ is a 2-spherical object and is generated by the 2-CY objects $\OO_{C_1},\OO_{C_2}(-1),\ldots,\OO_{C_n}(-1)$. 

It is a separate question if the subcategory $\Db_D(X)$ of objects supported on $D$ is a 2-CY category, \ie has Serre functor $[2]$. This holds if $D$ is of ADE type, see \autoref{cor:ade2cy}.
\end{remark}

\subsection{Curvelike decompositions}

\begin{definition} \label{def:minimal-divisor}
A $(-n)$-divisor $D$ is called \emph{minimal} 
if no $(-1)$-curves can be contracted from $D$, and no $(-2)$-curves can be twisted off $D$, \ie if
\begin{itemize}[topsep=0ex]
\item $D.C \neq 0$ for all $(-1)$-curves $C$ in $D$, and
\item $D.C \neq -1$ for all $(-2)$-curves $C$ in $D$.
\end{itemize}
\end{definition}

\begin{remark}
Any negative curve is a minimal $(-n)$-divisor, for $n>1$.
If a $(-n)$-divisor $D$ can be obtained from a $(-n)$-curve by blow-ups or spherical twists, then $D$ is said to be \emph{essentially a $(-n)$-curve}.
\end{remark}

\begin{lemma}
\label{lem:minimal-twistoff-comps}
A $(-n)$-divisor is minimal if and only if no $(-1)$-curve in $D$ can be contracted and no spherical component can be twisted off.
\end{lemma}

\begin{proof}
  The direction $(\Longleftarrow)$ is clear. For the converse, let $D$ be a $(-n)$-divisor having a spherical component $A$ but no contractible $(-1)$-curves. We have to show that $D$ is not a minimal $(-n)$-divisor, \ie there is a $(-2)$-curve that can be twisted off.
  
Assume the contrary, \ie $D.C \neq -1$ for all $(-2)$-curves $C \led D$. Then, as $D$ is $1$-connected, $(D-C).C \geq 2$ or, equivalently, $D.C \geq 0$.
In particular, by \autoref{cor:curvelike-2-spherical} the spherical component $A$ consists entirely of $(-2)$-curves, so we find $A.D \geq 0$. This is absurd as $A.D = A^2 + A.(D-A) = -1$.
\end{proof}

The above proof (with $D=A$) also works for the following statement; note that a spherical divisor consists of $(-2)$-curves only.

\begin{corollary} \label{cor:spherical-divisors-essential}
A spherical divisor is essentially a $(-2)$-curve.
\end{corollary}

\begin{definition} \label{def:curvelike-decomposition}
Let $D$ be a $(-n)$-divisor for some $n>1$. A \emph{curvelike decomposition} of $D$ is a sum $D = D_1 + \cdots + D_m$ of effective subdivisors 
such that all $D_i$ are $(-n_i)$-divisors with $-2 \geq -n_i \geq -n$ for all $i$ and $D_i.(D_{i+1}+\cdots+D_m)=1$ for  $i<m$.
\end{definition}

For a description of curvelike decompositions, we coin a term for the divisors obtained by successively blowing up a $(-n)$-curve:
\[ \intext{\tikzsimplechain{n}}, \qquad \text{called \emph{simple $(-n)$-chain}} . \]
We assume that such a simple $(-n)$-chain contains always the $(-n-1)$-curve, but not necessarily any $(-2)$-curve.

\begin{theorem}
\label{prop:curvelike-chopped}
Let $D$ be a $(-n)$-divisor on $X$ which is not the pullback of a $(-n)$-curve.
If $n \geq 2$ then there is a contraction $\pi\colon X\to Y$ to a smooth surface $Y$ with
$D = \pi^* D'$ and a 1-decomposition
  $D' = C_1 + \cdots + C_l + \cdots + C_m$
inducing a curvelike decomposition $D'=A+B$, $A = C_1+\cdots+C_l$ and either
\begin{itemize}
\item $A$ is a $(-k)$-curve with $-2 \geq -k \geq -n$ and $B$ is a $-(n+2-k)$-divisor,
\item or $A$ is a simple $(-n)$-chain with $C_1^2 = -n-1$ and $C_l^2=-1$ and $B$ is a $(-2)$-divisor.
\end{itemize}
\end{theorem}

\begin{remark}
The divisor of \autoref{ex:sph-3-2211} admits only curvelike decompositions $D=A+B$ with both $A$ and $B$ chains. 
So the second clause of the theorem above is necessary.
\end{remark}

This theorem yields the curvelike decomposition $D=A_1+\ldots+A_k$ of \autoref{thm:B} by applying it iteratively.

\begin{corollary}
A minimal $(-n)$-divisor which is not a curve has a non-trivial curvelike decomposition.
\end{corollary}

\begin{remark}
Let $D = \pi^*A+\pi^*B$ be a curvelike decomposition obtained from \autoref{prop:curvelike-chopped} with $A$ a simple $(-n)$-chain.
In particular, there is a $(-1)$-curve $C$ in $A$.
As $A$ is built from a $1$-decomposition, $C.D = C.A + C.B = 0+1=1$. In particular, $D$ does not come as a pullback of some divisor from the surface where $C$ is contracted.
\end{remark}

Before turning to the proof of \autoref{prop:curvelike-chopped}, we need the following lemma.

\begin{lemma} \label{lem:decomp-from-1decomposition}
Let $D$ be a negatively closed divisor with a $1$-decomposition $D=C_1+\cdots+C_m$. 
Suppose that for some $l<m$, all $A_i = C_1+\cdots+C_i$ are reduced trees, where $i=1,\ldots,l$. Then
\[
A_i^2+1 = A_i.D, \quad
A_i.(D-A_i)=1,   \quad\text{and}\quad
-1 \geq A_i^2 \geq D^2-1.
\]
\end{lemma}

\begin{proof}
We add the 1-decomposition pieces $1 = C_k.(D-C_1-\cdots-C_k)$ up to $i \leq l$:
\begin{align*}
  i &= \sum_{k=1}^i C_k.(D-C_1-\cdots-C_k)
     = \sum_k C_k.D - \sum_{k'<k} C_{k'}.C_k - \sum_k C_k^2 \\
    &= A_i.D - A_i^2 + \sum_{k'<k} C_{k'}.C_k
     = A_i.D - A_i^2 + i-1 ,
\end{align*}
where the final sum expands to $i-1$ because $A_i$ is a reduced tree with $i$ vertices. This implies the first formula.
The second one is obtained by plugging $D=(D-A_i)+A_i$ into the first one.
From this, we get
\[
0 > (D-A_i)^2 = D.(D-A_i)-A_i.(D-A_i) = D^2-D.A_i-1=D^2-A_i^2-2
\]
as $D$ is negatively closed. For the same reason also $A_i^2<0$ and, combining these two, the last statement follows.
\end{proof}

In \autoref{def:curvelike-decomposition} and \autoref{prop:curvelike-chopped}, we exclude $(-1)$-divisors, since these have no non-trivial curvelike decompositions:

\begin{proposition} \label{cor:exceptional-divisors}
A $(-1)$-divisor is essentially a $(-1)$-curve.
\end{proposition}

\begin{proof}
Let $D = C_1 + \cdots + C_m$ be a 1-decomposition of a $(-1)$-divisor.
Then by \autoref{lem:decomp-from-1decomposition}, $C_1$ can be either a $(-1)$-curve or a $(-2)$-curve. As $C_1.(D-C_1) = 1$, this curve can be blown down or twisted off, respectively.
\end{proof}

\begin{proof}[Proof of \autoref{prop:curvelike-chopped}]
Without loss of generality, we may assume that there are no $(-1)$-curves $C$ in $D$ that can be contracted, \ie with $D.C=0$.
Let $D = C_1+\cdots+C_m$ be a $1$-decomposition of a $(-n)$-divisor $D$ and set $A_i \coloneqq C_1+\cdots+C_i$ and $B_i \coloneqq D-A_i$, so that $C_i.B_i = 1$.
From \autoref{lem:decomp-from-1decomposition}, we know $-1 \geq C_1^2 = A_1^2 \geq D^2-1=-n-1$. There are three cases:
\begin{itemize}[itemsep=+0.5ex]
\item $C_1^2=-1$ is excluded, since it could be contracted (as $D.C_1=0$);
\item $-2 \geq C_1^2 \geq -n$, and $D= C_1 + (D-C_1)$ forms a curvelike decomposition: $D-C_1$ is well-connected and rigid by \autoref{lem:rigid-descends} and \ref{lem:1-decomposition-given}, so we are done;
\item $C_1^2=-n-1$, and then $A_1.B_1 = 1$ and $B_1^2=-1$.
\end{itemize}
So, we only have to deal with the last case, where a simple $(-n)$-chain will be constructed inductively.
By \autoref{rem:can-complete-1decomps}, we may assume in the following that there is no $(-2)$-curve $C$ in $D$ with $C.(D-C)=1$ as we could start the 1-decomposition of $D$ with such a curve (which could be twisted off $D$).

For the induction step, assume that for some $i$ we have:
\begin{itemize}[itemsep=1ex]
\item $C_1^2 = -n-1$ and $C_2^2=\ldots=C_i^2=-2$, so $A_i^2=-n-1$;
\item $B_i^2=-1$ and $A_i.B_i=1$.
\end{itemize}
In particular, we find that $i+1<m$: for $i+1=m$ we have $B_i=C_m$, and as $D.C_m=0$ we would get a contractible $(-1)$-curve $C_m$ in $D$. This we have excluded.

First note that by the assumption $D$ negatively closed we get
\begin{align*}
0 > (B_i-C_{i+1})^2 &= B_i.(B_i-C_{i+1})-C_{i+1}.(B_i-C_{i+1})  \\
                    &= B_i^2 - B_i.C_{i+1} - 1 = -2-B_i.C_{i+1} ,
\end{align*}
using also the $1$-decomposition and $B_i^2=-1$. So this implies $B_i.C_{i+1}\geq-1$.
Combining it again with $1 = C_{i+1}.(B_i-C_{i+1})$ from the $1$-decomposition, we get $C_{i+1}^2 \in \{-1,-2\}$.

Next, we will show $A_i.C_{i+1} \in \{0,1\}$.
To see this, note that by the 1-decomposition and the induction hypothesis
\begin{align*}
  \tag{$\ast$}
  \label{eq:next-intersection-AB}
A_{i+1}.B_{i+1}  & = (A_i+C_{i+1}).(B_i-C_{i+1}) \\
                 & =  A_i.B_i-A_i.C_{i+1} + C_{i+1}.B_{i+1} = 2-A_i.C_{i+1}.
\end{align*}
As $D$ is $1$-connected by \autoref{lem:curvelike-num1conn}, we get $1 \leq A_{i+1}.B_{i+1} = 2- A_i.C_{i+1}$, so $A_i.C_{i+1} \leq 1$.
On the other hand, as $A_i = (D - C_{i+1}) - B_{i+1}$, we get that
\[
A_i.C_{i+1} = (D-C_{i+1}).C_{i+1} - B_{i+1}.C_{i+1} \geq 1 - 1 = 0 ,
\]
since $D$ is $1$-connected.

\Case{1} $A_i.C_{i+1} = 0$.
Then we compute
\[
D.C_{i+1} = B_i.C_{i+1}=B_{i+1}.C_{i+1}+C_{i+1}^2=1+C_{i+1}^2 .
\]
If $C_{i+1}^2=-1$, then $C_{i+1}.D=0$, whereas if $C_{i+1}^2=-2$, then $C_{i+1}.(D-C_{i+1}) = 1$. By assumption $D$ is minimal, so both cases are impossible.

\Case{2} $A_i.C_{i+1} = 1$. 
Recall that $i+1<m$, so we distinguish the following two subcases.

\Subcase{2.1}  $C_{i+1}^2=-1$.
As we have here $A_{i+1}^2=-n$ and $A_{i+1}.B_{i+1} = 1$ from \eqref{eq:next-intersection-AB}, we get a curvelike decomposition $D=A_{i+1}+B_{i+1}$.

\Subcase{2.2} $C_{i+1}^2=-2$. 
Then $A_{i+1}^2=-n-1$, $A_{i+1}.B_{i+1}=1$ and therefore $B_{i+1}^2=-1$. As $A_{i+1}$ does not contain $(-1)$-curves either, we can start the whole argument with $i \mapsto i+1$.

\medskip

So eventually, we will end at \emph{Subcase 2.1} with $A_l=C_1+\cdots+C_l$ where
\begin{gather*} C_1^2 = -n-1, \quad
   C_2^2 = \cdots = C_{l-1}^2 = -2, \quad
   C_l^2 = -1, \quad \text{and} \\
(C_1 + \cdots + C_i).C_{i+1} = 1 \quad \text{for } 1 \leq i < l .
\end{gather*}
We will construct a subdivisor $A$ of $A_l$ which is the desired simple $(-n)$-chain. By the second property there is a unique $i_1$ such that $C_l.C_{i_1}=1$.
Likewise, as $C_{i_1}.(C_1+\cdots+C_{i_1-1})=1$ there is a minimal index $i_2$ such that $C_{i_1}.C_{i_2}=1$.
Proceed inductively until we reach $C_{i_{j+1}}=C_1$.
We claim that $A = C_1 + C_{i_j} + \cdots + C_{i_1} + C_l$ gives the desired decomposition $D= A+(D-A)$.

By the (minimal) choice of the indices, $A$ is indeed a reduced chain. Next we see that
$C_{i_j}.(D-C_1-C_{i_j}) = C_{i_j}.(C_2+\cdots+C_{i_j-1}) + C_{i_j}.(C_{i_j+1}+\cdots+C_m)$.
By minimality the first summand is 0, whereas the second summand is 1 by the $1$-decomposition.
The same holds for the remaining intersection numbers,
\ie for $j \geq k \geq 0$, 
setting $i_0 \coloneqq l$ and $i_{j+1} \coloneqq 1$
\[
\resizebox{\textwidth}{!}{%
  $\begin{split}
    C_{i_k}.\Bigl(D-\sum_{s=k}^{j+1} C_{i_s}\Bigr) 
& = C_{i_k}.(C_{i_k+1}+\cdots+C_m) 
  + \sum_{s=k}^{j} C_{i_k}.(C_{i_{s+1}+1}+\cdots+C_{i_{s}-1}) \\
& = 1+0 .
  \end{split}$
}
\]
By \autoref{rem:can-complete-1decomps}, $A$ comes from a $1$-decomposition of $D$. In particular, $A.(D-A)=1$ due to \autoref{lem:decomp-from-1decomposition}. Finally, by \autoref{lem:1-decomposition-given} and \ref{lem:rigid-descends}, $D-A$ is well-connected and rigid.
\end{proof}

\subsection{Cohomology of \boldmath{$(-n)$}-divisors}

The two subsequent subsections provide applications of \autoref{prop:curvelike-chopped}. Here, we justify our terminology again:
various cohomology groups of $(-n)$-divisors coincide with $(-n)$-curves.

\begin{proposition} \label{prop:Odual-hom-groups}
Let $D$ be a $(-n)$-divisor. Then
\begin{align*}
  \Hom^\bullet(\OO_D,\OO) \cong \kk^{n-1}[-2], &&& H^\bullet(\OO_D(K)) \cong \kk^{n-1}, \\
  \Hom^\bullet(\OO_D(D),\OO) \cong \kk[-1],   &&& H^\bullet(\OO_D(D+K)) \cong \kk[-1].
\end{align*}
\end{proposition}

\begin{proof}
The claims about $H^\bullet(\OO_D(K))$ and $H^\bullet(\OO_D(D+K))$ follow from their left-hand counterparts by Serre duality. Moreover, the statement about $\Hom^\bullet(\OO_D,\OO)$ follows from one for $\Hom^\bullet(\OO_D(D),\OO)$. To see this, apply $\Hom^\bullet(\OO_D(D),\blank)$ to $\OO \into \OO(D) \onto \OO_D(D)$ and get the triangle
\[
\Hom^\bullet(\OO_D(D),\OO) \to \Hom^\bullet(\OO_D,\OO) \to \Hom^\bullet(\OO_D,\OO_D) .
\]
By definition of $D$ to be a $(-n)$-divisor, $\Hom^\bullet(\OO_D,\OO_D) \cong \kk\oplus\kk^{n-1}[-2]$, hence
  $\Hom^\bullet(\OO_D(D),\OO) = \kk[-1] \iff \Hom^\bullet(\OO_D,\OO)=\kk^{n-1}[-2]$.

To show $\Hom^\bullet(\OO_D(D),\OO) \cong \kk[-1]$, we may assume that $D$ is not the pullback of a smaller $(-n)$-divisor.
By \autoref{prop:curvelike-chopped} there is a decomposition $D=A+B$ where $A$ is either a $(-m)$-curve with $-2 \geq m \geq n$ or $A = \intext{\tikzsimplechain{n}}$.
We claim that in the triangle 
\[
\Hom^\bullet(\OO_B(B),\OO) \from \Hom^\bullet(\OO_D(D),\OO) \from \Hom^\bullet(\OO_A(D),\OO)
\]
the right-hand term vanishes, so that the remaining isomorphism allows the reduction to a divisor $D'$ which is the pullback of a $(-n')$-curve $C$.
But then
\[ \Hom^\bullet(\OO_D(D),\OO) \cong \Hom^\bullet(\OO_{D'}(D'),\OO) \cong
   H^\bullet(\OO_C(C+K))^*[-2] \cong \kk[-1] ,
\]
using Serre duality and $C.(C+K) = -2$.
It remains to show the vanishing of $\Hom^\bullet(\OO_A(D),\OO)$ or of $H^\bullet(\OO_A(D+K))$ by Serre duality.
If $A$ is just a curve, then the vanishing follows from $A.(D+K)= A^2 + A.B -A^2-2 = -1$.

Otherwise, let $E=C_l$ be the single $(-1)$-curve in $A$ and $\pi\colon X \to Y$ the contraction of $E$.
Then $A.E=0$, and hence $A = \pi^* A'$ for some shorter chain $A' = \intext{\tikzsimplechain{n}}$.
Moreover, we have $B.E=1$, as it came from a $1$-decomposition of $D$, \ie $D.E = 1$. Therefore, $D = \pi^* D' - E$.
Finally, $K_X = \pi^* K_Y + E$.
Putting all this together, we find $\OO_A(D+K) = \pi^* \OO_{A'}(D'+K)$.

Let $E'$ be the $(-1)$-curve in $A'$ whose strict transform becomes the last $(-2)$-curve $C=C_{l-1}$ in $A$.
Again we compute that
\begin{align*}
E'.D' &= \pi^* E'.\pi^* D' = (C+E).(D+E) = C.D +1 \\
      &= C_l.(C_1+\cdots+C_{l-1}) + C_l^2 + C_l.(C_{l+1}+\cdots+C_m) + 1 = 1,
\end{align*}
as the first and the third summand are equal to $1$.
This allows us to proceed inductively, until we contract $A$ to a single $(-n)$-curve.
\end{proof}

\subsection{Simple chains as spherelike components}

The following result complements \autoref{prop:spherelike-decomposition-with-one-spherical}.

\begin{proposition} \label{prop:simple-chain-component}
Let $D=C_1+\cdots+C_m$ be a 1-decomposition of a $(-n)$-divisor such that $A = C_1+\cdots+C_l = \intext{\tikzchain{3,2,,2,1}}$ is a spherelike component for some $2\leq l<m$.
Then
\[ \Hom^\bullet(\OO_A(A-D),\OO_D) = \kk\oplus\kk[-1]\oplus\kk[-2] . \] 
\end{proposition}

In particular, given a spherelike, non-spherical component $A$ of $D$ as in the statement, the twist functor $\TTT_{\OO_A(A-D)}$ is not an autoequivalence, and $\TTT_{\OO_A(A-D)}(\OO_D) \neq \OO_{D-A}$. Roughly speaking, $A$ cannot be twisted off $D$.

Before going for the proof, we need a lemma extending \autoref{prop:spherelike-decomposition-with-one-spherical}(1).

\begin{lemma}
\label{lem:homOAB-OD}
Let $D$ be a $(-n)$-divisor with a spherelike component $A$, and let $B=D-A$.
Then
\[
 \Hom^\bullet(\OO_A(-B),\OO_D) \cong \kk \oplus H^\bullet(\OO_A(K-B))^*[-2]
        \cong     \begin{cases} \text{either } \kk \text{ or } \\
                                \kk \oplus \kk[-1] \oplus \kk[-2]. \end{cases}
\]
\end{lemma}

\begin{proof}
Applying $\Hom^\bullet(\OO_A(-B),\blank)$ to the ideal sheaf sequence of $D$ yields the triangle
\[
\Hom^\bullet(\OO_A(A),\OO) \to \Hom^\bullet(\OO_A(-B),\OO) \to \Hom^\bullet(\OO_A(-B),\OO_D) .
\]
By \autoref{prop:Odual-hom-groups}, the left term is isomorphic to $\kk[-1]$. Using Serre duality, the middle term is isomorphic to $H^\bullet(\OO_A(K-B))^*[-2]$.
Combining these facts yields the first claimed equivalence, as $\Hom(\OO_A(-B),\OO_D)\neq0$.

For the second isomorphism of the statement, apply $\Hom^\bullet(\blank,\OO_D)$ to the decomposition sequence for $A+B$:
\[
\Hom^\bullet(\OO_A(-B),\OO_D) \from \Hom^\bullet(\OO_D,\OO_D) \from \Hom^\bullet(\OO_B,\OO_D) :
\]
Therefore $\kk \cong \Ext^2(\OO_D,\OO_D) \onto \Ext^2(\OO_A(-B),\OO_D)$. This completes the proof, because $\chi(\OO_A(K-B)) =\chi(\OO_A) + A.K - A.B = 1 + 0 -1 = 0$.
\end{proof}

\begin{proof}[Proof of \autoref{prop:simple-chain-component}]
According to \autoref{lem:homOAB-OD}, it suffices to show $H^0(\OO_A(K-B)) \neq 0$ instead.
Fix a 1-decomposition $D = C_1+\cdots+C_m$, and denote $D_i \coloneqq C_1+\cdots+C_i$. Then $A=D_l$.
Because we assume $D=A+B$ to be as in \autoref{prop:curvelike-chopped}, we have $\OO_A(K-B) = \OO_{D_l}(K-D+D_l)$.
We now proceed inductively on $D_i$ until we reach $i=l$.

\step{Step $i=1$:} 
Since $C_1$ is a $(-3)$-curve, we have $C_1.K=1$, and by the $1$-decomposition $C_1.(D-C_1)=1$:
\[
H^\bullet(\OO_{C_1}(K-D+C_1)) \cong H^\bullet(\OO_{C_1}) \cong \kk.
\]

\step{Step $1 < i < l$:}
We use the decomposition sequence for $D_{i-1} + C_i$, tensored with $\OO(K-D+D_i)$:
\[
   H^\bullet(\OO_{D_{i-1}}(K-D+D_{i-1})) \to 
   H^\bullet(\OO_{D_i}(K-D+D_i)) \to
   H^\bullet(\OO_{C_i}(K-D+D_i)) .
\]
By the induction hypothesis, the left term is isomorphic to $\kk$. The right term vanishes, since $C_i.(C_{i+1}+\cdots+C_m)=1$ from the $1$-decomposition and $C_i.K=0$ as $C_i$ is a $(-2)$-curve. Hence, the middle term is isomorphic to $\kk$.

\step{Step $i=l$:}
We look at the same triangle of cohomology spaces as in the previous step for $i=l$.
The left term is still isomorphic to $\kk$. But $C_l.K=-1$ as $C_l$ is a $(-1)$-curve, so the right term is isomorphic to $\OO_{\IP^1}(-2) \cong \kk[-1]$.
From this follows that the middle term is
\[
H^\bullet(\OO_A(K-B)) = H^\bullet(\OO_{D_l}(K-D+D_l)) = \kk \oplus \kk[-1]
\]
as claimed.
\end{proof}

\section{Spherelike divisors}
\label{sec:spherelike-divisors}

\noindent
Exceptional divisors, \ie $(-1)$-divisors, can always be dismantled to $(-1)$-curves by way of contractions and spherical twists; see \autoref{cor:exceptional-divisors}. In this section, we look at spherelike divisors, \ie $(-2)$-divisors.

We start with a simple observation: an easy way for checking that a divisor $D$ is spherelike is via the numerical criterion \autoref{cor:numerical-criterion}, specifying a 1-decomposition and a negative filtration, and testing $D.K_X=0$. And by \autoref{lem:self-squares}, $D.K_X=0$ if and only if the average of self-intersection numbers of all curve components is $-2$.

We recall central notions from \cite[\S 4]{HKP}. Let $D$ be a spherelike divisor and $D'$ any other divisor. Then by Serre duality 
\[
\Hom(\OO_D(D'),\OO_D(D'+K_X)) \cong \Ext^2(\OO_D,\OO_D)^* \cong \kk.
\]
Therefore, we have a non-zero map, unique up to scalars, which can be completed to the \emph{asphericity triangle}
\[
\OO_D(D') \to \OO_D(D'+K_X) \to Q_{\OO_D(D')} ,
\]
whose last term $Q_{\OO_D(D')}$ we call the \emph{asphericity} of $\OO_D(D')$.

\begin{definition}
Let $D$ be a spherelike divisor and $D'$ an arbitrary divisor.
The \emph{spherical subcategory} of the spherelike object $\OO_D(D')$ is
\[
\Db(X)_{\OO_D(D')} \coloneqq \lorth Q_{\OO_D(D')} =
          \{ M \in \Db(X) \mid \Hom^\bullet(M,Q_{\OO_D(D')})=0 \} .
\]
\end{definition}

\begin{proposition}[{\cite[Thm.~4.4 \& 4.6]{HKP}}]
Let $D$ and $D'$ be divisors on $X$, with $D$ spherelike. Then $\smash{\Db(X)_{\OO_D(D')}}$ is the unique maximal full triangulated subcategory of $\Db(X)$ containing $\OO_D(D')$ as a spherical object.
\end{proposition}

\subsection{Spherical subcategories and blow-ups}
First we have a look how the spherical subcategory changes under blow-ups.
Blowing up a surface $\pi\colon X'\to X$ in a single point, with exceptional divisor $E$, there are the following semi-orthogonal decompositions, see for example \cite[Prop. 11.8]{Huybrechts}:
\[
\Db(X') = \sod{\OO_E(E),\pi^*\Db(X)} = \sod{\pi^*\Db(X),\OO_E}.
\]

\begin{proposition}
\label{prop:asphericities-blowup}
Let $\pi \colon X' \to X$ be the blow-up of a surface $X$ in a point $P$ with exceptional divisor $E$.
For a spherelike divisor $D$ on $X$, the asphericities of $D$ and $D' = \pi^*D$ fit into the triangle
\[
\pi^* Q_{\OO_D} \to Q_{\OO_{D'}} \to \OO_{D'} \otimes \OO_{E}(E)
\]
satisfying a dichotomy for the derived tensor product $T = \OO_{D'} \otimes \OO_{E}(E)$:

\begin{tabular}{@{}l@{$\implies$}l@{}}
  $P \notin \supp(D)$ & $T = 0$ and $\Db(X')_{\OO_{D'}} = \sod{\OO_E(E),\pi^*\Db(X)_{\OO_D}}$; \\
  $P \in \supp(D)$ & $T = \OO_E(E) \oplus \OO_E(E)[1]$ and $\pi^* \Db(X)_{\OO_D} \subseteq \Db(X')_{\OO_{D'}}$.
\end{tabular}
\end{proposition}

\begin{proof}
We compare the pullback of the asphericity triangle of $\OO_D$ with that of $\OO_{D'}$ in the following diagram, 
where we set $K = K_X$, $K'=K_{X'}$, $Q = Q_{\OO_D}$ and $Q' = Q_{\OO_{D'}}$ and use $\pi^*K=K'-E$ and $E.K' = E^2$:
\begin{equation}
\tag{$Q$}
\label{eq:connecting-qs}
\begin{tikzcd}[column sep=0.5em, row sep=1.0ex]
\OO_{D'} \ar[rr] \ar[rd] & \ar[d, phantom, "\displaystyle \hypertarget{eq:cone}{(\ast)}"] 
& \OO_{D'}(K'-E) \ar[rrr] \ar[ld] &&& \pi^* Q \ar[llldd, dashrightarrow, in=20, out=-130]  \\
& \OO_{D'}(K')  \ar[ld] \ar[rd] \\
\OO_{D'} \otimes \OO_E(E) && Q' \ar[ll, dashrightarrow, in=-10, out=200]
\end{tikzcd}
\end{equation}
We now show that (\hyperlink{eq:cone}{$\ast$}) commutes, hence induces the claimed triangle by the octahedral axiom.
For this, note first that 
\begin{equation}
\tag{$D'$}
\label{eq:D'}
\OO_{D'}(K'-E) \to \OO_{D'}(K') \to \OO_{D'}\otimes\OO_E(E)
\end{equation}
is the ideal sheaf sequence of $E$ tensored with $\OO_{D'}(K')$.
Tensoring the ideal sheaf sequence of $D'$ with $\OO_E(E)$ and using $E.D'=0$, we get the triangle 
\[
\OO_E(E) \xto{\alpha} \OO_E(E) \to \OO_{D'}\otimes\OO_E(E) .
\]
Depending on the position of $P$, $\alpha$ is either an isomorphism or zero, hence
\[
\OO_{D'}\otimes\OO_E(E) \cong
\begin{cases}
\OO_E(E) \oplus \OO_E(E)[1] & P \in \supp(D);\\
0                             & P \not\in \supp(D).
\end{cases}
\]
In both cases, applying $\Hom(\OO_{D'},\blank)$ to \eqref{eq:D'} yields an isomorphism
\[
\Hom^\bullet(\OO_{D'},\OO_{D'}(K'-E)) \isom \Hom^\bullet(\OO_{D'},\OO_{D'}(K')),
\]
as $\Hom^\bullet(\OO_{D'},\OO_E(E)) = \Hom^\bullet(\OO_D,\pi_*\OO_E(E)) = 0$. Hence (\hyperlink{eq:cone}{$\ast$}) commutes.

We are left to show $\pi^* \Db(X)_{\OO_D} \subset \Db(X')_{\OO_{D'}}$ if $P \in \supp(D)$.
For this, let $A \in \Db(X)_{\OO_D}$, \ie $\Hom^\bullet(A,Q) = 0$.
Then applying $\Hom^\bullet(\pi^*A,\blank)$ to the triangle connecting the asphericities in \eqref{eq:connecting-qs}, we find
\[
\Hom^\bullet(A,Q) \to \Hom^\bullet(\pi^*A,Q') \to \Hom^\bullet(\pi^*A,\OO_E(E) \oplus \OO_E(E)[1]) . 
\]
The left term vanishes by assumption and the right term does because of $\pi_*\OO_E(E) = 0$.
Hence $\Hom^\bullet(\pi^*A,Q')=0$, so $\pi^*A \in \Db(X')_{\OO_{D'}}$.
\end{proof}

\begin{example}[{\cf \cite[Ex.~5.6]{HKP}}]
\label{ex:single-blowup}
Let $Y$ be a surface containing a $(-2)$-curve $C$ and $\psi\colon X \to Y$ be the blow-up in a point $P$.
Then the spherical subcategory of $\psi^* \OO_C$ is
\[
\Db(X)_{\psi^* \OO_C} = 
\begin{cases}
 \psi^* \Db(Y) & P \in C;\\
\Db(X) & P \not\in C.
\end{cases}
\]
\end{example}

The following proposition shows that the spherical subcategory also keeps track of more complicated blow-up situations.
We only treat the next case after \autoref{ex:single-blowup}. A general statement seems possible, but would have to take care of the combinatorial structure of iterated blow-ups.

\begin{proposition}
\label{prop:iterated-blowup}
Let $\pi\colon X' \to X$ be the blow-up of a surface $X$ in a point $P$ with exceptional divisor $E'$.
Let $D$ be a divisor on $X$ of type \intext{\tikzchain{3,1}}, and let $D' = \pi^* D$.
Finally, let $Y$ be the surface obtained by contracting the $(-1)$-curve in $D$ and $\pi' \colon X' \to X \to Y$.
Then one of the following three cases occurs, each distinguished by the  spherical subcategory of $\OO_{D'}$:
\begin{enumerate}
\item \label{it:Poutside}
  $D'$ is of type \tikzchain{3,1}\ \ 
  if $P\notin D$;
  then $\Db(X')_{\OO_{D'}} = \psi'^*\Db(Y')$ where $\psi'\colon X' \to Y'$ is the contraction of the $(-1)$-curve in $D'$;
\item \label{it:Pon3}
  $D'$ is of type \tikzchain{1,4,1}\ \ 
  if $P$ lies on the $(-3)$-curve, but not on the $(-1)$-curve;
  then $\Db(X')_{\OO_{D'}} = \pi'^*\Db(Y)$; 
\item \label{it:Pon1}
  $D'$ is of type \tikzchain{3,2,1}\ or \tikzchainmult{4,1,2}{1,2,1}\ \ 
  if $P$ lies on the $(-1)$-curve;
  then $\Db(X')_{\OO_{D'}} = \sod{\OO_C(-1)} \oplus \pi'^*\Db(Y)$ where $C$ is the $(-2)$-curve in $D'$.
\end{enumerate}
\end{proposition}

Let $M \in \Db(X')$.
In the statement and proof of \autoref{prop:iterated-blowup}, we denote by $\sod{M}$ the smallest triangulated full subcategory of $\Db(X')$ which contains $M$ and is closed under taking direct summands.

\begin{proof}
Note that $X$ is the blow-up $\psi\colon X \to Y$ of a surface $Y$ in a point $Q$ lying on a $(-2)$-curve with exceptional divisor $E$ and total transform $D$; this is the setting of \autoref{ex:single-blowup}.
Hence Case \eqref{it:Poutside} becomes a direct application of \autoref{ex:single-blowup}: 
Let $\psi' \colon X' \to Y'$ be the contraction of the $(-1)$-curve in $D'$.
Then we can write $\OO_{D'}$ as $(\psi')^* \OO_C$ where $C$ is a $(-2)$-curve on $Y'$,
and therefore $\Db(X')_{\OO_{D'}} = \psi'^* \Db(Y').$
 
In the other cases, \autoref{prop:asphericities-blowup} produces the following triangle of asphericities:
\begin{equation}
\tag{Q}
\label{eq:qs}
\pi^*\OO_E(E) \oplus \pi^*\OO_E(E)[1] \to Q_{\OO_{D'}} \to \OO_{E'}(E') \oplus \OO_{E'}(E')[1] .
\end{equation}
In Case \eqref{it:Pon3}, the composition $\pi'\colon X'\xxto{\pi} X \xxto{\psi} Y$ is obtained from blowing up $Q$ and $\psi(P)$ in either order and, in particular, $\Hom^\bullet(\OO_{E'}(E'),\pi^*\OO_E(E))=0$ by $\supp(E')\cap\supp(\pi^*E)=\varnothing$.
Hence $Q_{\OO_{D'}}$ is the direct sum of the outer terms of Triangle \eqref{eq:qs} and
  $\Db(X')_{\OO_{D'}} = \lorth(\pi^*\OO_E(E) \oplus \OO_{E'}(E')) = \pi'{}^*\Db(Y)$,
using 
  $\Db(X') = \sod{\OO_{E'}(E'), \pi^*\OO_E(E), \pi'^*\Db(Y)}$.

In Case \eqref{it:Pon1}, we have a closer look at the degree increasing morphism in Triangle \eqref{eq:qs}. One can check that $\Hom^\bullet(\OO_{E'}(E'),\pi^*\OO_E(E)) \cong \kk \oplus \kk[-1]$, so the following arrows are possible:
\[
\begin{tikzcd}[row sep=1ex, ampersand replacement=\&]
\OO_{E'}(E')[1] \ar[d, phantom, "\oplus"] \ar[r, "e_2"] \ar[dr, "\alpha"] \& \pi^* \OO_E(E)[2] \ar[d, phantom, "\oplus"]\\
\OO_{E'}(E')  \ar[r, "e_1"'] \& \pi^* \OO_E(E)[1] 
\end{tikzcd}
\]
Consider the curve $C \coloneqq \pi^*E-E'$. Then $C^2=-2, C.E'=1,C.\pi^*E = -1$. 

\Step{1} $\OO_C(\pi^*E) \in \Db(X')_{\OO_{D'}} = \lorth Q_{\OO_{D'}}$.

\noindent
Note $\OO_{C}(\pi^*E) \in \sod{\OO_{E'}(E'),\pi^*\OO_E(E)}$ from a decomposition sequence.
Now $\Hom^\bullet(\OO_C(\pi^*E),\OO_{D'}) = 0$ implies
  $\Hom^\bullet(\OO_C(\pi^*E),Q_{\OO_{D'}}) = 0$
by \cite[Thm.\ 4.7]{HKP}, and thus it suffices to show the former vanishing.
Consider the decomposition sequence
$
 \OO_C(-E') \into \OO_{\pi^*E} \onto \OO_{E'}
$
and observe $C.(-E') = C.\pi^*E$, so that $\OO_C(-E') \cong \OO_C(\pi^*E)$.
Applying $\Hom^\bullet(\blank,\OO_{D'})$ gives
\[
\Hom^\bullet(\OO_C(\pi^*E),\OO_{D'}) 
\from
\Hom^\bullet(\OO_{\pi^*E},\OO_{D'})
\from
\Hom^\bullet(\OO_{E'},\OO_{D'}) .
\]
Using $\Db(X') = \sod{\pi'{}^*\Db(Y),\pi^*\OO_E,\OO_{E'}}$ or by direct calculation,  the middle and right terms are zero. Hence the left term vanishes too, as claimed.
Thus the spherical subcategory of $\OO_{D'}$ contains the spherical object $\OO_C(\pi^*E)$.

\Step{2} $e_1 \neq0$ and $e_2\neq0$.

\noindent
Assume the contrary. Then the asphericities have as direct summand $\OO_{E'}(E')$ if $e_1=0$, or $\pi^* \OO_E(E)[1]$ if $e_2=0$. 
Now observe that
\begin{itemize}
\item $\Hom^\bullet(\OO_C(\pi^*E),\OO_{E'}(E')) \neq 0$ because of the extension corresponding to the decomposition sequence $\OO_{E'}(E') \into \pi^*\OO_E(E) \onto \OO_C(\pi^*E)$;
\item $\Hom^\bullet(\OO_C(\pi^*E),\pi^*\OO_E(E)) \cong \Hom^\bullet(\OO_C,\OO_{E'+C}) \neq 0$ because we find $\chi(\OO_C,\OO_{E'+C}) = -C.(E'+C)=1 \neq 0$.
\end{itemize}
In particular, we would get that $\OO_C(\pi^*E) \not\in \lorth Q_{\OO_{D'}}$, contradicting Step 1.

\Step{3} $\Db(X')_{\OO_{D'}} = \bigsod{\sod{\OO_{E'}(E'),\pi^*\OO_E(E)} \cap \lorth \OO_F(F), \pi'{}^*\Db(Y)}$.

\noindent
This will follow from an explicit formula for $Q_{\OO_{D'}}$.
Consider the divisor $F = \pi^*E + E' = C + 2E'$ of type \smash{\tikzchainmult{2,1}{1,2}}.
Tensoring the decomposition sequence of $F$ with $\OO(F)$ yields
$
\pi^* \OO_E(E) \into \OO_F(F) \onto \OO_{E'}(E').
$
As the arrows $e_i$ are non-zero, one can check that $Q_{\OO_{D'}}$ is quasi-isomorphic to 
\[
Q_{\OO_{D'}} \cong
\begin{cases}
Q_s \coloneqq \OO_F(F) \oplus \OO_F(F)[1] & \text{if $\alpha=0$,}\\
Q_n \coloneqq \cone(\OO_F(F) \xto{t} \OO_F(F)) & \text{if $\alpha\neq0$,}
\end{cases}
\]
with $t \in \Hom^\bullet(\OO_F(F),\OO_F(F)) = \Hom^\bullet(\OO_F,\OO_F) \cong \kk[t]/t^2$ and $\deg(t)=0$.

Even though these two possibilities differ, their left orthogonals inside $\sod{\OO_{E'}(E'),\pi^*\OO_E(E)}$ coincide.
For the inclusion $\lorth Q_s \subseteq \lorth Q_n$,
let $A\in \lorth Q_s$, \ie $\Hom^\bullet(A,\OO_F(F))=0$. This forces $\Hom^\bullet(A,Q_n) = 0$ by applying $\Hom^\bullet(A,\blank)$ to the triangle $\OO_F(F) \to \OO_F(F) \to Q_n$.

On the other hand, let $B\in \lorth Q_n$ and suppose $0 \neq \beta \in \Hom^\bullet(B,\OO_F(F))$.
Applying $\Hom^\bullet(B,\blank)$ to the triangle defining $Q_n$, we get an isomorphism $t^* \colon \Hom^\bullet(B,\OO_F(F)) \to \Hom^\bullet(B,\OO_F(F))$.
Now $t^*(\beta) = t \beta \neq 0$, but $t^*(t\beta) = t^2 \beta = 0$, a contradiction.

\Step{4} 
We prove $\sod{\OO_{E'}(E'),\pi^*\OO_E(E)} \cap \lorth \OO_F(F) = \sod{\OO_C(-1)}$ via tilting.

\noindent
Consider the triangulated category $\CC = \sod{\EE_0,\EE_1} \coloneqq \sod{\OO_{E'}(E'),\pi^*\OO_E(E)}$.
The generators satisfy the following properties: (a) $\EE_0,\EE_1$ are exceptional; (b) $\Hom^\bullet(\EE_1,\EE_0)=0$; (c) $\Hom^\bullet(\EE_0,\EE_1)=\kk\oplus\kk[-1]$; (d) non-zero morphisms $\EE_0\to\EE_1$ are injective. The results of \cite{HP1} apply to exceptional sequences of this kind. Specifically, \cite[Prop.~1.7]{HP1} states that $\CC \cong \Db(\End(T)\hh\mod)$, where $T$ is the iterated universal extension of the exceptional sequence.

$\OO_F(F)$ is the unique non-trivial extension of $\pi^*\OO_E(E)$ by $\OO_{E'}(E')$, hence the tilting object is $T = \pi^*\OO_E(E) \oplus \OO_F(F)$ and $A \coloneqq \End(T)$ is the quiver algebra

{\centering
\begin{tikzpicture}
\node (A) at (0,0) {$1$};
\node (B) at (3,0) {$2$};
\draw [->, bend left=15] (A) to node[above] {$\scriptstyle\alpha$} (B);
\draw [->, bend left=15] (B) to node[below] {$\scriptstyle\beta$} (A);
\clip (-1,-0.7) rectangle (4,0.7);
\end{tikzpicture}
\par}

\vspace{-1ex}
\noindent
with relation $\beta\alpha=0$.
Under the equivalence $\CC \isom \Db(A\hh\mod)$, 
the $\pi^*\OO_E(E)$ and $\OO_F(F)$ become the projective modules $Ae_1$ and $Ae_2$, respectively, where $e_i$ is the idempotent to the vertex $i$. 
Moreover, $\OO_C(-1)$ is sent to the simple module $S_1$ associated to the vertex $1$.
One can check that for a finite-dimensional quiver algebra $B$ and an idempotent $e_i$ corresponding to a vertex $i$ holds
\[
Be_i^\perp = 
\sod{S_j \mid j \neq i}.
\]
In particular, $Ae_2^\perp = \sod{S_1}$.
Going back to $\CC$, this implies that $\OO_F(F)^\perp = \sod{\OO_C(-1)}$.
As $\OO_C(-1)$ is spherical and in particular a Calabi--Yau object, we 
conclude that $\lorth \OO_F(F) = \sod{\OO_C(-1)}$ as well.
By similar reasoning, $\OO_C(-1)$ is two-sided orthogonal to $\pi'{}^* \Db(Y)$,
completing the proof.
\end{proof}

\begin{remark}
The algebra $A$ of Step 4 of the proof appears in two well-known series: it is the derived-discrete algebra $\Lambda(1,2,0)$, see \cite{Vossieck,BGS}, and it is the Auslander algebra of $\kk[t]/t^2$, see \cite{HP1, HP2}. 
\end{remark}

\begin{remark}
One can check that in Case \eqref{it:Pon1}, the asphericities distinguish between \tikzchain{3,2,1} and \tikzchainmult{4,1,2}{1,2,1}, namely,
\[
Q_{\OO_{D'}} =
\begin{cases}
\OO_F(F) \oplus \OO_F(F)[1] & \text{if $D' = $ \tikzchainmult{4,1,2}{1,2,1};}\\
\cone(\OO_F(F) \xto{t} \OO_F(F)) & \text{if $D' = $ \tikzchain{3,2,1}.}
\end{cases}
\]

Moreover, in Case \eqref{it:Pon1} we have $\Db(X')_{\OO_{D'}} = \sod{\OO_C(-1)} \oplus \pi'^*\Db(Y)$, an orthogonal decomposition.
Also, $\sod{\OO_C(-1)}$ is neither weakly admissible in $\Db(X')$, i.e.\ its inclusion does not admit any adjoint, nor is it equivalent to $\Db(Z)$ for any smooth projective surface $Z$.
As a consequence, the same applies to $\Db(X')_{\OO_{D'}}$.
\end{remark}

\subsection{Spherical subcategories and decompositions}

If $D$ is a spherelike divisor with a curvelike decomposition $D=A+B$ (\autoref{def:curvelike-decomposition}), then $A$ and $B$ are spherelike divisors themselves, and we emphasise this by calling $D=A+B$ a \emph{spherelike decomposition}.

The remaining part of this section is dedicated to the behaviour of the spherical subcategory under spherelike decompositions.
\autoref{prop:spherelike-decomposition-with-one-spherical} and \cite[Lem.~2.2]{HKP2} together imply the following lemma.

\begin{lemma}
\label{lem:twisting-geometric}
Let $D$ be a spherelike divisor.
If $A \led D$ is a spherical component, then $\Db(X)_{\OO_D} = \TTinv_{\OO_A(-B)}(\Db(X)_{\OO_B})$ where $B=D-A$.
\end{lemma}

\begin{proposition} \label{prop:asphericities}
Let $D=A+B$ be a spherelike decomposition. Then the asphericity of $\OO_D$ occurs in a triangle $R_A \to Q_{\OO_D} \to R_B$ where
\begin{align*}
  R_A &= \begin{cases}
            \mathrlap{Q_{\OO_A(-B)}, } \hphantom{mmmmmmmmmmmmm}
                                          & \text{if } H^0(\OO_B(D)) = 0 \\
            \OO_A(K-B)\oplus\OO_A(-B)[1], & \text{if } H^0(\OO_B(D)) = \kk
         \end{cases} \\
  R_B &= \begin{cases}
           \mathrlap{Q_{\OO_B}, }     \hphantom{mmmmmmmmmmmmm}
                                          & \text{if } H^0(\OO_A(K-B)) = 0 \\
            \OO_B(K)\oplus\OO_B[1],       & \text{if } H^0(\OO_A(K-B)) = \kk
         \end{cases}
\end{align*}
\end{proposition}

\begin{proof}
Note that the above clauses for $R_A$ and $R_B$ exhaust all cases. For $R_B$ this follows from \autoref{lem:homOAB-OD}.
For $\OO_B(D)$ one can see this by the triangle $H^\bullet(\OO_A(A)) \to H^\bullet(\OO_D(D)) \to H^\bullet(\OO_B(D))$ from the decomposition $D=A+B$ and using that $H^\bullet(\OO_A(A)) \cong \kk[-1]$ and $H^0(\OO_D(D))=0$.
We recall that $Q_{\OO_D}$ is computed as the cone of the canonical map $\OO_D \xxto{\omega} \OO_D(K)$.
Using the decomposition sequence for $D=A+B$ and its $\OO(K)$-twist, we get two triangles linked by $\omega$:
\begin{equation} \label{asphericity-diagram} \tag{$\ast$}
\begin{gathered}  
\begin{tikzcd}
  \OO_A(-B) \ar[r, "\iota"] \ar[d, dashrightarrow, "\alpha"] &
  \OO_D     \ar[r, "\pi"] \ar[d, "\omega"] &
  \OO_{B}    \ar[d, dashrightarrow, "\beta"] \\
  \OO_A(K-B) \ar[r, "\iota'"]     & \OO_D(K) \ar[r, "\pi'"]   & \OO_{B}(K)
\end{tikzcd}
\end{gathered}
\end{equation}
We already know $\Hom(\OO_A(-B),\OO_B(K))=\Ext^2(\OO_B,\OO_A(-B))^*=0$ from \autoref{lem:spherelike-decomp}(2) and Serre duality; in particular $\pi'\omega\iota=0$ and this implies that $\omega$ extends to a map of triangles. In fact, $\Ext^{-1}(\OO_A(-B),\OO_B(K))=0$ ensures that $\omega$ determines the resulting morphisms $\alpha$ and $\beta$ uniquely. Taking cones, $R_A=\cone(\alpha)$ and $R_B=\cone(\beta)$, we get a triangle $R_A \to Q_{\OO_D} \to R_B$. However, there are various cases, depending on whether $\alpha\neq0$ or $\beta\neq0$.

If $\alpha\neq0$, then this morphism is a multiple of $\omega_{\OO_A(-B)}$, and hence its cone is $R_A = Q_{\OO_A(-B)} = Q_{\OO_A}(-B)$. On the other hand, if $\alpha=0$, then the triangle is split and $R_A$ is the direct sum given in the statement of the proposition. The same reasoning applies to $\beta$ and $R_B$.
We now look at $R_A$:
\begin{align}
  \alpha=0 &\iff \iota'\alpha=\omega\iota=0 \\
           &\iff \Hom(\OO_B,\OO_D(K)) \neq 0 \\
           &\iff H^1(\OO_B(D)) \neq 0        \\
           &\iff H^0(\OO_B(D)) = \kk
\end{align}

\iffstep{1}
$\iota'$ is injective, and the left-hand square of \eqref{asphericity-diagram} commutes.

\iffstep{2}
apply $\Hom(\blank,\OO_D(K))$ to the top triangle of \eqref{asphericity-diagram} to get
\[
 0 \to \Hom(\OO_B,\OO_D(K)) \to \Hom(\OO_D,\OO_D(K)) \xxto{\iota^*} \Hom(\OO_A(-B),\OO_D(K)) .
\]
If $\Hom(\OO_B,\OO_D(K)) = 0$, then $\iota^*$ is injective, mapping $\omega \mapsto \omega \iota = \iota'\alpha \neq 0$. On the other hand, if $\Hom(\OO_B,\OO_D(K)) \neq 0$, then $\iota^*=0$ because we have $\Hom(\OO_D,\OO_D(K)) = \Ext^2(\OO_D,\OO_D)^* = \kk$.

\iffstep{3}
applying $\Hom(\blank,\OO_B)$ to $\OO(-D) \into \OO \onto \OO_D$ yields the following snippet of the long exact sequence, proving $H^1(\OO_B(D)) \cong \Hom(\OO_B,\OO_D(K))$:
\[
\resizebox{\textwidth}{!}{ 
\begin{tikzcd}[cramped, column sep = 2em, row sep = 2.4ex, ampersand replacement=\&]
 \Ext^1(\OO,\OO_B) \ar[d, equal] \ar[r] \& \Ext^1(\OO(-D),\OO_B) \ar[r] \ar[d, equal] \& \Ext^2(\OO_D,\OO_B) \ar[d, equal] \ar[r] \& \Ext^2(\OO,\OO_B) \ar[d, equal] \\
  0 = H^1(\OO_B)              \ar[r] \& H^1(\OO_B(D))          \ar[r]            \& \Hom(\OO_B,\OO_D(K))^*         \ar[r] \& H^2(\OO_B) = 0
\end{tikzcd}
}
\]

\iffstep{4}
By Riemann--Roch, $\chi(\OO_B(D)) = -\frac{1}{2}(B^2+B.K)+B.D = 1-1 = 0$, so
 $H^0(\OO_B(D)) \cong H^1(\OO_B(D)) \cong \Hom(\OO_B,\OO_D(K))$.

Now we look at the other cone $R_B$:
\begin{align}
  \beta=0 &\iff \beta\pi = \pi'\omega = 0     \tag{1'} \\
          &\iff \Hom(\OO_D,\OO_A(K-B)) \neq 0 \tag{2'} \\
          &\iff H^0(\OO_A(K-B)) = \kk         \tag{3'}
\end{align}

\iffstep{1'}
$\pi$ is surjective, and the right-hand square of \eqref{asphericity-diagram} commutes.

\iffstep{2'}
apply $\Hom(\OO_D,\blank)$ to the bottom triangle of \eqref{asphericity-diagram} and get
\[ 0 \to \Hom(\OO_D,\OO_A(K-B)) \to \Hom(\OO_D,\OO_D(K)) \xxto{\pi'_*} \Hom(\OO_D,\OO_B(K)) . \]
$\pi'_*$ is injective $\iff \Hom(\OO_D,\OO_A(K-B))=0$, as $\Hom(\OO_D,\OO_D(K)) = \kk$.

\iffstep{3'} this is Serre duality applied to \autoref{lem:homOAB-OD}.
\end{proof}

\begin{remark}
Let $D$ be spherelike and consider its asphericity triangle
 $\OO_D \to \OO_D(K) \to Q$.
Taking cohomology and combining it with \autoref{prop:Odual-hom-groups} yields $H^\bullet(Q) = \Hom^\bullet(\OO,Q) = 0$, \ie $\OO_X\in\Db(X)_{\OO_D}$.
\end{remark}

\begin{example} \label{ex:sph-3-2211}
Let $D = 2B+C+C'+E+E' = \tikztreemult{2,2,1,1}{3}{1,1,1,1}{2}$ such that $B^2=-3, C^2=C'^2=-2, E^2=E'^2=-1$ and $D$ is rational.

$D$ is spherelike by \autoref{cor:numerical-criterion}: a negative filtration is $B,C,C',B,E,E'$, a 1-decomposition is $B,C,C',E,E',B$ (only crucial that $B$ comes first and last or next-to-last), and finally $D.K=0$ is immediate.


The algorithm proving \autoref{prop:curvelike-chopped} produces a spherelike decomposition out of a given 1-decomposition. The 1-decomposition $B,C,C',E,E',B$ from above yields
  $D = (B + C + C' + E) + (B + E')$.
Starting with the 1-decomposition $B,C,E,C',E',B$, we get
  $D = (B + C + E) + (B + C' + E')$.

Because of $E.D = E'.D = 1 \neq 0$, neither $(-1)$-curve can be contracted to yield a smaller spherelike divisor; see \autoref{sub:modify-blowup}. Similarly, $C.D = C'.D = 0 \neq 1$ means that neither of the $(-2)$-curves can be twisted off $D$; see \autoref{sub:modify-twisting}. Therefore, $D$ is a minimally spherelike divisor.

About the asphericity of $D$:
we employ the criterion of \autoref{prop:asphericities} with the  spherelike decomposition $D = A + A'$ with $A = B + C + C' + E$ and $A' = B + E'$. So we have to compute $H^0(\OO_{A'}(D))$ and $H^0(\OO_A(K-A'))$. For the former, the $\OO(D)$-twisted decomposition sequence of $A' = B + E'$
\[ 
\begin{tikzcd}[row sep=2ex]
\OO_{E'}(D-B) \ar[r, hook] \ar[d, equal] & \OO_{A'}(D) \ar[r, two heads] & \OO_B(D) \ar[d, equal]\\
\OO_{E'}                        &                   & \OO_B(-2)
\end{tikzcd}
\]
yields $H^0(\OO_{A'}(D))=\kk$. Thus the map $\alpha\colon \OO_A(-A') \to \OO_A(K-A')$ is zero, and hence $R_A = \OO_A(-A')[1] \oplus \OO_A(K-A')$.
In order to calculate $H^0(\OO_A(K-A'))$, we use the decomposition sequence for $A = (C+C'+E) + B$ and note that $C+C'+E$ is a disjoint union:
\[ 
\begin{tikzcd}
\OO_C(-1) \oplus \OO_{C'}(-1) \oplus \OO_E(-1) \ar[r, hook] & \OO_D \ar[r, two heads] & \OO_B  ,
\end{tikzcd}
\]
and twist it by $\OO(K-A')$:
\[ 
\begin{tikzcd}
\OO_C(-2) \oplus \OO_{C'}(-2) \oplus \OO_E \ar[r, hook] & \OO_D(K-A') \ar[r, two heads] & \OO_B ,
\end{tikzcd}
\]
using $K.C=0, K.E=-1, K.B=1$ and $A'.C = A'.E =  A'.B=1$. We find $0\neq H^0(\OO_D(K-A'))$,
and therefore $\beta\colon \OO_{A'}\to\OO_{A'}(K)$ is also zero, forcing $R_{A'} = \OO_{A'}[1] \oplus \OO_{A'}(K)$. The asphericity of $\OO_D$ thus sits in the triangle
\[ 
\begin{tikzcd}
\OO_A(-A')[1] \oplus \OO_A(K-A') \ar[r] &  Q_{\OO_D} \ar[r] &  \OO_{A'}[1] \oplus \OO_{A'}(K) . 
\end{tikzcd}
\]
We take the cohomology exact sequence of this triangle:
\[
\begin{tikzcd}
0 \ar[r] &  \OO_A(-A')  \ar[r] & h^{-1}(Q_D) \ar[r] \ar[d, phantom, "\bar\omega"{name=M}] & \OO_{A'}  \ar[r, no head, out=-20, in=0, to=M] \\
         &  \OO_A(K-A') \ar[r] \ar[out=180, in=170, from=M] &  h^0(Q_D)   \ar[r] & \OO_{A'}(K) \ar[r] & 0.
\end{tikzcd}
\]
The map $\bar\omega$ really is induced from $\omega\colon\OO_D\to\OO_D(K)$, again by $\alpha=0$ and $\beta=0$. The only common component of $A$ and $A'$ is the $(-3)$-curve $B$, so $\omega$ must be nontrivial there. We have $B.(K-A') = B.K - B.(B+E') = 1 - (-3+1) = 3$. Now $\hom(\OO_B,\OO_B(3))=4$, but $\omega$ has to vanish on the intersections of $B$ with the other curves in $A$ ($C.(K-A')=C'.(K-A')=-2, E.(K-A')=0$, so no poles allowed). Prescribing these three zeroes, there is a unique map $\OO_B\to\OO_B(3)$. Splitting the above long exact sequence into short exact sequences, we get:
\[ 
\begin{tikzcd}[row sep=0ex]
  (a) &   \OO_A(-A')        \ar[r, hook] &  h^{-1}(Q_D)  \ar[r, two heads] &  \OO_{E'}(-B)       \\
  (b) &   \OO_{E'}(-B)       \ar[r, hook] &  \OO_{A'}    \ar[r, two heads] &  \OO_B              \\
  (x) &   \OO_B(-C-C'-E)    \ar[r, hook] &  \OO_A       \ar[r, two heads] &  \OO_{C+C'+E}        \\
  (c) &   \OO_B             \ar[r, hook] &  \OO_A(K-A') \ar[r, two heads] &  \OO_{C+C'+E}(K-A')  \\
  (d) &   \OO_{C+C'+E}(K-A') \ar[r, hook] &  h^0(Q_D)    \ar[r, two heads] &  \OO_{A'}(K)        
\end{tikzcd}
\]
Here (a)--(b)--(c)--(d) splice to give the long exact sequence. (b) and (x) are decomposition sequences. (c) is (x) twisted by $K-A'$, using $B.(A'-K) = -3$ for $\OO_B(-C-C'-E) = \OO_B(-3) = \OO_B(A'-K)$. The last term of (c) is
\[ \OO_{C+E+E'}(K-A') = (\OO_C\oplus\OO_E\oplus\OO_{E'})(K-A') = \OO_C(-1)\oplus\OO_E(-2)\oplus\OO_{E'}(-2) . \]
The two further twisted decomposition sequences
\[
\begin{tikzcd}[row sep=0ex]
   \OO_A(-A')        \ar[r, hook] &  \OO_{B+C+C'+E+E'}(-B)  \ar[r, two heads] &  \OO_{E'}(-B)    \\
   \OO_{C+C'+E}(K-A') \ar[r, hook] &  \OO_{B+C+C'+E+E'}(K)    \ar[r, two heads] &  \OO_{A'}(K) 
\end{tikzcd}
\]
show
\begin{align*}
  h^{-1}(Q_D) &= \OO_{B + C + C' + E + E'}(-B) , \\
  h^0(Q_D)   &= \OO_{B + C + C' + E + E'}(K) .
\end{align*}
The degrees of these line bundles on $B+C+C'+E+E'$ differ.
At this point, it seems hard to compute the spherical subcategory
$\Db(X)_{\OO_D} = \orth Q_D$ explicitly. We do know $\OO_D,\OO_X\in\Db(X)_{\OO_D}$ and $\lorth h^{-1}(Q_D) \cap \lorth h^0(Q_D) \subseteq \Db(X)_{\OO_D}$.
\end{example}

\section{Negative definite divisors and rational singularities}
\label{sec:negative-definite}

\noindent
We recall a few facts about surface singularities, as can be found in \cite{Artin,Artin2,BHPV,Badescu}.
For modern proofs of the contraction results, see \cite[\S A.7 and \S 4.15f]{Reid}.

\begin{definition} 
A normal surface $Y$ is called a \emph{rational singularity} if
there is a resolution of singularities $\pi\colon X \to Y$ such that $\pi_* \OO_X = \OO_Y$, \ie $\pi$ has connected and acyclic fibres ($R^i\pi_*\OO_X=0$ for $i>0$).
\end{definition}


\begin{proposition}[{\cite[Thm.~2.3]{Artin}}] \label{prop:Artin-contraction}
Let $D$ be a reduced, connected divisor on an algebraic surface $X$.
Then $D$ can be contracted to a point $P$ on an algebraic surface $Y$ with $\chi(Y) = \chi(X)$ if and only if $D$ is negative definite, and all effective divisors supported on $D$ are Jacobi rigid.

If this holds, then $P$ is a rational singularity.
\end{proposition}

\begin{remark} \label{rem:Artin-contraction}
The original phrasing of the second condition had $\chi(\OO_{D'})\geq1$ for all $D'$ supported on $D$, instead of $H^1(\OO_{D'})=0$ (Jacobi rigid). Moreover, $H^0(\OO_{D'}(D'))$ vanishes for such divisors, too. To see this, note that such $D'$ is negative definite as well, in particular can be negatively filtered, so $H^0(\OO_{D'}(D'))=0$ by \autoref{lem:negatively-filtered-implies-subscheme-rigid}.
Hence any effective divisor supported on a configuration yielding a rational singularity is automatically rigid.

There is a related result stating that $D$ is contractible in the analytic category if and only if $D$ is negative definite, see \cite{Grauert}.
\end{remark}

\begin{definition} \label{def:Znum}
Let $\pi\colon X \to Y$ be a resolution of the normal surface singularity $P\in Y$. So the exceptional locus $\pi^{-1}(P) = \cup_i C_i$ is a union of projective curves.
The \emph{numerical (fundamental) cycle} is the minimal divisor $\Znum = \sum z_i C_i$ such that $z_i > 0$ and $\Znum.C_i \leq 0$ for all $i$.
\end{definition}

\begin{remark} \label{rem:nonbirational-property}
The condition defining numerical cycles makes sense for arbitrary divisors: $D$ is called \emph{anti-nef} if $C.D \leq 0$ for all curves $C \leqd D$. This property is obviously numerical in the sense of \autoref{def:divisor-properties}. One can check that the property is birational, because $D$ allows contraction of a $(-1)$-curve $E\led D$ if and only if $E.D=0$.

Now consider the stronger variant $C.D < 0$ for all curves $C\leqd D$. It is also numerical but certainly not birational: a reduced $(-1)$-curve satisfies the condition, but the blow-up \intext{\tikzchain{1,2}} does not. We remark that this property forces $D$ to be negative definite: let $M$ the intersection matrix of $D\red$; then the intersection matrix of $D$ is $TMT$ where $T$ is the diagonal matrix of curve coefficients of $D$. By $C.D < 0$, all column sums of $TMT$ are negative, \ie $-TMT$ is a symmetric, strictly diagonally dominant matrix with positive diagonal entries, hence positive definite \cite[Cor.~1.22]{Varga}.
\end{remark}

\begin{remark}
It is possible that the dual graph of a resolution of a normal surface singularity is not a tree. A singularity is called \emph{arborescent} if some resolution produces a tree (this property then holds for any good resolution); see \cite[\S 4]{Arborescent}. Rational singularities have this property, see \eg \autoref{prop:Znum-curvelike}.
\end{remark}

\begin{definition}
A rational singularity $Y$ is called an \emph{ADE singularity} if there exists a crepant resolution of singularities $\pi\colon X \to Y$, \ie $K_X = \pi^* K_Y$.
\end{definition}

This definition is anachronistic, but the most convenient one for us. As is well-known, there are many characterisations of these singularities, leading to a lot of equivalent terminology, such as rational double point, simple surface singularity, and they are often named after du Val or Klein; see \cite{Durfee}.

\begin{proposition} \label{prop:Znum-curvelike}
If $X \to Y$ is a resolution of a rational singularity, then $\Znum$ is $(-n)$-divisor for some $n>0$.
In particular, $\Znum$ is spherelike if and only if $Y$ is an ADE singularity.
\end{proposition}

\begin{proof}
As noted in \autoref{rem:Artin-contraction}, any effective divisor supported on the exceptional locus is rigid. Especially this holds for $\Znum$.
Moreover by \cite[Prop.~4.12]{Reid}, $\chi(\OO_{\Znum}) = 1$, so $\Znum$ is well-connected.

If $Y$ is an ADE singularity, then $\Znum^2 = -2$, \ie $\Znum$ is spherelike.
The converse implication holds by \cite[Thm.~3.31]{Badescu}.
\end{proof}

\begin{corollary}
\label{cor:ade2cy}
Let $X \to Y$ be a crepant resolution of an ADE singularity. Then $\Znum$ is spherical. Moreover, the subcategory $\Db_{\Znum}(X)$ of objects set-theoretically supported on $\Znum$ is a 2-Calabi--Yau category.
\end{corollary}

\begin{proof}
We already know that $\Znum$ is spherelike.
The singular surface $Y$ is Gorenstein and the resolution is crepant, therefore there is an open subset $U\subset X$ containing the exceptional locus such that $\omega_X|_U \cong \OO_U$. In particular, $\omega_X |_{\Znum} \cong \OO_{\Znum}$. This shows $\Znum$ is spherical, and also that $M\otimes\omega_X \cong M$ for any $M$ supported on $\Znum$. Hence $\Db_{\Znum}(X)$ has Serre functor $-\otimes\omega_X[2]=[2]$, \ie is a 2-Calabi--Yau category.
\end{proof}

\begin{example} \label{ex:nonreduced-spherical} \label{ex:d4}
\tikztreemult{2,2,2}{2}{1,1,1}{2} is the numerical cycle of a minimal resolution of a $D_4$-singularity, it is a non-reduced spherical divisor. See \cite[p.~96]{BHPV} for the complete list of the numerical cycles for the ADE singularities.
\end{example}

\begin{proposition}
\label{prop:Znum-max-ndiv}
Let $D$ be a negative definite, $(-n)$-divisor that can be contracted.
Then $D \leqd \Znum$. In particular if $D$ contracts to a rational singularity, then $\Znum$ is the maximal well-connected and rigid divisor with support $\supp D$.
\end{proposition}

The proof uses that the numerical cycle can be computed by Laufer's algorithm, and one can see a $1$-decomposition as a special case of it.

\begin{proof}
By \cite[Prop.~4.1]{Laufer}, the numerical cycle can be computed recursively: Start with $Z_0 \coloneqq C$ for some curve $C \led D$. Given $Z_i$, compute $Z_i.C'$ for all curves $C' \led D$. If $Z_i.C' > 0$ for some $C'$, set $Z_{i+1} \coloneqq Z_i+C'$; else $Z_i.C' \leq 0$ for all $C'$, and then $\Znum = Z_i$.

Using a $1$-decomposition of $D = C_1+\cdots+C_m$ backwards, this algorithm yields $Z_i = C_{m-i}+\cdots+C_m$. In particular, $D = Z_{m-1} \leqd \Znum$.
\end{proof}

\begin{example} \label{ex:ADE-graphs}
Let $T$ be a tree of $(-2)$-curves. As a reduced divisor, $T$ is spherical by \autoref{cor:spherical-divisors-essential}.
If $T$ forms an ADE graph, then there is a unique maximal spherical divisor on $T$, the numerical cycle $\Znum$ of \autoref{def:Znum}.

The reverse implication holds true as well: if $T$ is not an ADE graph, there is no maximal spherical divisor. For example, consider the following two spherical divisors on a $\tilde D_4$-configuration of $(-2)$-curves:
\[
D_1 = \tikztreemult{2,2,2}{2}{1,1,1}{2}
\quad\text{and}\quad
D_2 = \tikztree{2,2,2,2}{2}
\]
The smallest divisor $D$ with $D_1\leqd D$ and $D_2\leqd D$ has $D^2=0$.
In particular, any divisor containing both $D_1$ and $D_2$ cannot be rigid, as it contains $D$ as a subdivisor.
\end{example}

\begin{proposition}
\label{prop:negative-definite-spherelike}
A negative definite, spherelike divisor can be contracted to either a smooth point or an ADE singularity and in the latter case is the pullback of a spherical divisor.
\end{proposition}

\begin{proof}
If the spherelike divisor $D$ contains no $(-1)$-curves, then it has to be a configuration of $(-2)$-curves, since the average among self-intersection numbers of all curves in $D$ is $-2$ by \autoref{lem:self-squares}. By \autoref{cor:curvelike-2-spherical} $D$ is spherical. Now $D$ is negative definite, and hence an ADE configuration of $(-2)$-curves, see for example \cite[\S 3]{Durfee}.
It is well-known that such a configuration can be contracted to an ADE singularity, which in turn gives $D \leqd \Znum$ by \autoref{prop:Znum-max-ndiv}.

If $D$ is not spherical, it contains a $(-1)$-curve $E$. As $D$ is $1$-connected, $(D-E).E \geq 1$, so $D.E \geq 0$.
On the other hand $D$ is negative definite, so $-1 \geq (D+E)^2 = -2 + 2D.E -1$ and hence $D.E \leq 1$.

If $D.E=0$, contract $E$ and start over with a smaller spherelike divisor.

If $D.E=1$, then $(D+E)^2=-1$. By negative definiteness, $D+E$ is negatively filtered, and so $H^0(\OO_{D+E}(D+E)) = 0$ by \autoref{lem:negatively-filtered-implies-subscheme-rigid}. On the other hand, the decomposition sequence yields a triangle
$
H^\bullet(\OO_E(-D)) \to H^\bullet(\OO_{D+E}) \to H^\bullet(\OO_D).
$
As $D.E=1$, we get $H^\bullet(\OO_{D+E})=\kk$.
Altogether, we find that $D+E$ is a $(-1)$-divisor and, moreover, can be contracted to a smaller $(-1)$-divisor $D'$, since $(D+E).E=0$. 
Starting this proof again with this $D'$ shows that there has to be a $(-1)$-curve which can be contracted, inductively yielding a smooth point.

Otherwise, we have $D.E=0$ throughout, so that $D$ eventually becomes the pullback of a spherical divisor.
\end{proof}

\begin{remark}
With this proof one can also show that a negative definite $(-1)$-divisor can be iteratively contracted to a smooth point.
\end{remark}

We end this section with two more examples, first of a divisor which contracts to an elliptic singularity, and then a spherelike divisor which is not negative definite.

\begin{example}[{\cite[Ex.~4.20]{Nemethi}}]
\label{ex:elliptic-sing}
Consider the minimal resolution of the surface singularity $\{x^3+y^3+z^4=0\} \subset \kk^3$, which is a minimally elliptic singularity, in particular not rational. Its numerical cycle $\Znum$ is
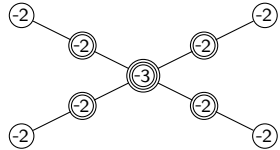
\begin{center}
\tikz
{
\node[curve,double] (m)  at ( 0  , 0)        {}; 
\node[curve]        (d1) at (-1.6,-0.8)      {\printnum{2}};
\node[curve,double] (d2) at (-0.8,-0.4)      {\printnum{2}};
\node[curve,double] (d3) at ( 0.8,-0.4)      {\printnum{2}};
\node[curve]        (d4) at ( 1.6,-0.8)      {\printnum{2}};
\node[curve]        (u1) at (-1.6, 0.8)      {\printnum{2}};
\node[curve,double] (u2) at (-0.8, 0.4)      {\printnum{2}};
\node[curve,double] (u3) at ( 0.8, 0.4)      {\printnum{2}};
\node[curve]        (u4) at ( 1.6, 0.8)      {\printnum{2}};
\draw (d1) -- (d2) -- (m) -- (d3) -- (d4);
\draw (u1) -- (u2) -- (m) -- (u3) -- (u4);
\node[tripl]              at ( 0  , 0)        {}; 
\node[curve,double] (m2)  at ( 0  , 0)        {\printnum{3}}; 
}
\end{center}
The reduced divisor $D = (\Znum)_{\red}$ is a negative definite $(-3)$-divisor which can be twisted off to a single $(-3)$-curve. By contrast, an easy computation shows $\chi(\OO_{\Znum})=0$.
Note that $\Znum$ is not $1$-decomposable, as $(\Znum-C).C=2$ for all curves $C \led \Znum$. Moreover, $\Znum$ is a negative definite divisor that is not Jacobi rigid.
\end{example}

\begin{example} \label{ex:spherelike-not-pullback}
Let $D = D_n \coloneqq B + E + C_1 + \cdots + C_n = \tikztree{3,2,,2}{1}$ where $B^2=-3, E^2=-1$ and $C_i^2=-2$. Then $D.K=0$ is obvious and $D$ is reduced, so pruning leaves yields both a 1-decomposition and a negative filtration.
Alternatively, one can first twist off the $(-2)$-curves, and then blow down the remaining $(-1)$-curve. In particular, $D$ is essentially a $(-2)$-curve.

We remark that $D_n$ for $n>1$ is not contractible to a rational singularity. This holds even in the analytic category, since $(B+3E+C_1+C_2)^2 = 2$ shows that $D_n$ is not negative definite.
The divisor $D_1$ is not the pull-back of any divisor, but contracts to a smooth point.
\end{example}

\section{Classification of minimal $(-n)$-divisors}
\label{sec:classification}

\subsection{Graph-theoretical algorithm}
$(-n)$-divisors have a discrete, or combinatorial, flavour. More precisely, fixing the self-intersection number $-n$ and the topological type as a graph $T$, there are only finitely many building blocks, \ie minimal $(-n)$-divisors. Here, we deal with exhausting these graphs algorithmically. There remains the question which of those graphs actually occur as dual intersection graphs of effective divisors (which are then necessarily $(-n)$-divisors), and this is taken up in \ref{sub:graphs-to-divisors}.

Formally speaking, we consider weighted graphs with multiplicities below. Nonetheless, we will speak of `curves' instead of `vertices', and `self-intersection number' instead of `weight'.

\begin{enumerate}
\item Start with all curves of multiplicity 1 and unknown self-intersection.
\item Form a list of partially defined divisors on $T$ (some curves may not yet have assigned a self-intersection number), increasing in each step one multiplicity $kC \mapsto (k+1)C$, such that the 1-decomposability condition is met. If $k=1$ this fixes $C^2$. 
\item For the remaining entries, fill unassigned self-intersection numbers in all possible ways, admitting a negative filtration and satisfying $D.K=n-2$ or, equivalently, $D^2=-n$.
\item Remove divisors having a $(-1)$-curve which can be contracted or a $(-2)$-curve which can be twisted off.
\end{enumerate}

Graphs surviving the final step possess a 1-decomposition and a negative filtration, and have self-intersection number $-n$. Hence, if they are dual graphs of divisors, these are $(-n)$-divisors by \autoref{cor:numerical-criterion}.
The resulting list then needs to be condensed, because it will contain multiple incarnations of the same divisor.
Moreover, after Step (2) the resulting list will be infinite in general. Still any $(-n)$-divisor on a given graph will eventually be covered by this algorithm.
The algorithm becomes more efficient if the following intermediate checks are also carried out during and after Step (2):
\begin{enumerate}
\item[(2')] Remove an entry from that list if there is a subdivisor violating the negatively closed property, \eg \intext{\tikzchain{1,1}} or \intext{\tikzchain{1,2,1}}.
\item[(2'')] Remove divisors having a $(-1)$-curve which can be contracted or a $(-2)$-curve which can be twisted off.
\end{enumerate}

Since by \autoref{cor:exceptional-divisors}, $(-1)$-divisors can always be worked down to $(-1)$-curves, here we investigate spherelike divisors.

\begin{example}
If $T$ is a tree with two, three or four vertices, then no spherelike divisor on $T$ is minimal. We show this if $T$ is the four chain; the reasoning for the other cases is similar. The list of Step (2), cleaned up by (2'), is:
\begin{center}
   \tikzchainmult{?,1,?,?}{1,2,1,1}, \hfill
   \tikzchainmult{?,1,2,?}{1,2,2,1}, \hfill
   \tikzchainmult{?,1,2,?}{1,3,2,1} 
\end{center}
\noindent
together with the reduced chain \intext{\tikzchain{?,?,?,?}} which, by \autoref{prop:reduced-spherelike}, cannot be a minimally spherelike divisor. The $(-1)$-curves in the first and third divisors can be contracted. The $(-2)$-curve in the second one can be twisted off.
\end{example}

\begin{example} \label{ex:minimally-spherelike-5}
We list all minimally spherelike divisors on five curves. Their existence as divisors (not just graphs) follows from \autoref{prop:existence-of-divisors}, or can also easily checked by hand.

\begin{center}
  \hfill \tikztypeD{2}{1}{3,1,3}{2,2,1}
  \hfill \tikztreemult{2,2,1,1}{3}{1,1,1,1}{2}
  \hfill \tikztypeDmirror{2}{2}{3,1,2}{2,2,1}
  \hfill{ }

\medskip
  \hfill \tikzchainmult{3,1,3,1,3}{1,2,2,2,1}
  \hfill \tikzchainmult{2,3,1,2,3}{1,2,3,2,1}
  \hfill{ }
\end{center}
\end{example}

\begin{proposition} \label{prop:reduced-spherelike}
A reduced spherelike divisor is essentially a $(-2)$-curve.
\end{proposition}

\begin{proof}
We will show that $D$ has a leaf (\ie a curve component intersecting only one other curve) of self-intersection $-1$ or $-2$: such a curve can be blown down or twisted off, obtaining a smaller divisor with the same properties.

For a contradiction, assume that $D$ is a reduced negatively closed tree such that $C^2\leq -3$ for all leaves $C$. Put $L$ for the subdivisor consisting of all, say $l$, leaves, and let $I \coloneqq D-L$ be the complement of all inner curves.

Note that $I\neq0$ unless $D$ is either a single curve or two curves intersecting transversally in a point, neither of which is possible under the assumption. Now $I\neq0$ implies four facts: First,  $I.L=l$, as $D$ is reduced and each leaf intersects exactly one inner curve, with multiplicity 1. Second, $L^2\leq-3l$, as $L$ is a disjoint union of $l$ curves $C$ with $C^2\leq-3$. Third, $I^2<0$, as $D$ is negatively closed. Fourth, $l\geq2$.
We obtain a contradiction from $D^2=-2$:
\begin{align*}
 D^2 &= (I+L)^2 
      = I^2 + 2I.L + L^2
      =    I^2 + 2l + L^2 \\
    &\leq I^2 + 2l -3l
      =    I^2 - l
      <    -l \qedhere
\end{align*}
\end{proof}

\begin{remark}
The proof of the proposition shows a bit more: if $D$ is a reduced $(-3)$-divisor that is not a chain, then it is essentially a $(-3)$-curve. The provision is necessary, an example is \intext{\tikzchain{3,1,3}}.
\end{remark}

\subsection{From graphs to divisors} \label{sub:graphs-to-divisors}
The previous subsection produces a list of the weighted graphs which can occur as the dual intersection graphs of divisors with prescribed properties. However, it is a subtle problem to decide which graphs can actually be realised by divisors. 

\begin{example} \label{ex:violates-HIT}
The following graph cannot be realised on any surface:
\[
\tikzchain{2,2,2,2,1,3,3,1,2,2,2,2}
\]
Suppose the contrary and let $D$ be a rational reduced divisor with the dual intersection graph above.
One can easily check that $D$ is $1$-decomposable (as it is reduced) and negatively filtered. Moreover, $D^2 = -2$, so $D$ would be a spherelike divisor.

But if we iteratively blow down the $(-1)$-curves next to the two middle curves, then we will end up after five steps with $C_1+C_2=\tikzchainpositive{2,2}$.
Such a configuration cannot exist on any surface, since then $C_1.(C_1-2C_2)=0$ but $(C_1-2C_2)^2>0$ contradicting the Hodge Index Theorem \cite[Thm.~V.1.9]{Hartshorne}.
\end{example}

Many graphs can be realised, however. For the next statement, we consider a weighted tree $T$, with vertices denoted by $C$, and weights $C^2$. Write $v(C)$ for the valency of a vertex $C$, \ie the number of vertices adjacent to $C$. We introduce a local quantity $\sigma(C)$ for each vertex, which measures excess positivity, and a global quantity $b(T)$ counting the number of bad vertices:
\begin{align*}
  \sigma(C) &\coloneqq C^2 + v(C) \\
  b(T)      &\coloneqq \#\{ C \mid C^2 > -v(C) \} = \#\{ C \mid \sigma(C) \geq 1 \} .
\end{align*}
Moreover, the \emph{distance} $d(C,C')$ of two vertices is the number of edges in the shortest path connecting $C$ and $C'$. 

\begin{proposition} \label{prop:existence-of-divisors}
Let $T$ be a finite weighted tree with badness $b = b(T)$, and bad vertices $C_1,\ldots,C_b$.
There exists a reduced divisor on a rational surface with dual intersection graph $T$ if $b\leq 1$ or one of the following hold:
\begin{center}
\begin{tabular}{@{}ll@{ }l@{}}%
(1) & $b = 2$, & $d(C_1,C_2)=1$, and $\sigma(C_1)=1$ or $\sigma(C_2)=1$ or \\
    &&  \phantom{$d(C_1,C_2)=1$, and} $\sigma(C_1) = \sigma(C_2) = 2$; \\
(2) & $b = 2$, & $d(C_1,C_2)=2$, and $\sigma(C_1) = \sigma(C_2) = 1$; \\
(3) & $b = 3$, & and $C_1,C_2,C_3$ form a 3-chain with $\sigma(C_1) = \sigma(C_2) = \sigma(C_3) = 1$.
\end{tabular}
\end{center}
\end{proposition}

\begin{remark}
Minimal examples for the case (1) of the proposition are
   \intext{\tikzchainpositive{0,\text{$\ast$}}} and \intext{\tikzchainpositive{1,1}}.
For the cases (2) and (3):
   \intext{\tikzchainpositive{0,\tinyminus2,0}} and
   \intext{\tikzchainpositive{0,\tinyminus1,0}}.
Note that of these two, the former can be obtained from the latter by blowing up. Moreover, that last chain can be obtained by blowing up the intersection point of \intext{\tikzchainpositive{1,1}}, ignoring multiplicities.
\end{remark}

\begin{proof}
It suffices to realise those trees as divisors $D$ which satisfy the equality $C^2 = -v(C)$, \ie $\sigma(C)=0$, for all but the bad vertices: any other weighted tree with smaller prescribed self-intersection numbers can be obtained by blowing up appropriate interior points on the curves of $D$.

Assume $b=0$. Then all leaves have weight $-1$. Thus we can contract each leaf in the numerical sense, \ie remove it and increase the weight of its neighbour by $1$. 
Moreover, the condition also guarantees that this process can be iterated, stopping at a single vertex of weight $0$. We can revert this process on any surface with a $0$-curve, \eg $\IP^1\times\IP^1$.

Assume $b=1$. We apply the same procedure, but now we end up with the single bad vertex of weight $m=\sigma(C_1)>0$.
Let $F_m \coloneqq \IP(\OO_{\IP^1}\oplus\OO_{\IP^1}(m))$ be the Hirzebruch surface containing a smooth rational curve $L$ with $L^2=m$. Again the process can be reverted, starting with $L$ and ending with a tree $D$ of rational curves whose dual intersection graph is $T$.

Assume $b=2$. We can numerically contract all vertices except for the two bad vertices $C_1,C_2$ and the path between them, obtaining
  \intext{\tikzchainpositive{$x$,\tinyminus2,,\tinyminus2,$y$}}
with $x \coloneqq \sigma(C_1)-1$ and $y\coloneqq \sigma(C_2)-1$.
If $C_1$ and $C_2$ are adjacent, \ie $d(C_1,C_2)=1$, then each of the three cases can be realised on some Hirzebruch surface: \intext{\tikzchainpositive{0,$m$}} on $F_m$ and \intext{\tikzchainpositive{1,1}} on $\IP^2$.
The chain \intext{\tikzchainpositive{0,\tinyminus2,0}} is a double blow-up of \intext{\tikzchainpositive{1,1}}.

Assume $b=3$. After contracting we arrive at \intext{\tikzchainpositive{0,\tinyminus1,0}} already mentioned in the remark above.
\end{proof}

\begin{question}
Which graphs can be realised as dual graphs of divisors?

  In a related vein: given a rational rigid (or negatively closed) divisor $D$ on some surface, can it be realised on a rational surface?
\end{question}

\addtocontents{toc}{\protect\setcounter{tocdepth}{-1}}   
\section*{Acknowledgements}
\addtocontents{toc}{\protect\setcounter{tocdepth}{1}}     

We want to thank all anonymous referees for many valuable comments. Their reports have improved this work considerably. We also want to thank Klaus Hulek, Martin Kalck and Malte Wandel.

\input{biblio-surfaces}

\bigskip
\noindent
\emph{Contact:} \texttt{andreas.hochenegger@sns.it, david.ploog@uis.no}

\end{document}

%% file: header-surfaces.tex
\usepackage{amsmath,amssymb,amsthm,amsfonts,bbm}
\usepackage{etex}
\usepackage{mathtools}

\usepackage[dvispsnames]{xcolor}
\definecolor{addcolor}{RGB}{0,60,60}
\definecolor{changecolor}{RGB}{120,0,0}

\usepackage{tikz}
\usetikzlibrary{shapes} 
\usetikzlibrary{calc}   
\usetikzlibrary{cd}

\usepackage{ifthen}
\usepackage{mdwlist}
\usepackage{ulem}
\usepackage{enumitem}
\setitemize[0]{leftmargin=*}
\usepackage{booktabs}

\usepackage{hyperref}
\hypersetup{
    colorlinks,
    linkcolor={red!50!black},
    citecolor={blue!50!black},
    urlcolor={blue!80!black}
}
\usepackage{aliascnt} 

\usepackage{setspace} 

\newcommand{\ie}{i.e.\ }
\newcommand{\eg}{e.g.\ }
\newcommand{\cf}{cf.\ }

\newcommand{\fixedwidthtabular}{ \noindent \begin{tabular}{@{} p{0.08\textwidth} @{} p{0.92\textwidth} @{} } }

\newcommand{\bib}[6]{{\bibitem{#2} #3, {\emph{#4},} #5#6.}}

\newcommand{\TTT}{\mathsf{T}\!}  

\newcommand{\TTinv}{\TTT^{\:-1}}

\newcommand{\Db}{\mathcal{D}^b}  
\newcommand{\CC}{\mathcal{C}}
\newcommand{\DD}{\mathcal{D}}
\newcommand{\EE}{\mathcal{E}}
\newcommand{\FF}{\mathcal{F}}

\newcommand{\OO}{\mathcal{O}}

\newcommand{\IP}{\mathbb{P}}

\newcommand{\kk}{\mathbbm{k}} 

\renewcommand{\mod}{\mathsf{mod}}

\newcommand{\hh}{\text{-}}  

\newcommand{\sod}[1]{{\langle #1 \rangle}}

\newcommand{\bigsod}[1]{{\big\langle #1 \big\rangle}} 

\newcommand{\Znum}{Z_{\mathrm{num}}} 
\newcommand{\leqd}{\preceq}
\newcommand{\led}{\prec}

\newcommand{\orth}{^\perp}

\newcommand{\lorth}{{}^\perp}

\newcommand{\blank}{\:\cdot\:} 
\newcommand{\red}{_{\text{red}}}

\DeclareMathOperator{\coker}{coker}

\DeclareMathOperator{\Hom}{Hom}
\DeclareMathOperator{\End}{End}

\DeclareMathOperator{\Ext}{Ext}
\DeclareMathOperator{\ext}{ext}

\DeclareMathOperator{\id}{id}

\DeclareMathOperator{\Pic}{Pic}

\DeclareMathOperator{\cone}{cone}

\DeclareMathOperator{\supp}{supp}

\newcommand{\cdtext}[1]{ \text{#1} }
\newcommand{\cdstext}[2]{{ \begin{array}{@{}c@{}} \text{#1} \\ \text{#2} \end{array} }}
\newcommand{\num}{^{n}}
\newcommand{\closed}{^{c}}
\newcommand{\closednum}{^{{cn}}}

\newcommand{\isom}{ \text{{\hspace{0.48em}\raisebox{0.8ex}{${\scriptscriptstyle\sim}$}}}
                    \hspace{-0.65em}{\rightarrow}\hspace{0.3em}} 
\newcommand{\embed}{\hookrightarrow}
\newcommand{\into}{\hookrightarrow}
\newcommand{\onto}{\twoheadrightarrow}
\newcommand{\from}{\leftarrow}

\newcommand{\xto}[1]{\xrightarrow{#1}}
\newcommand{\xxto}[1]{\xrightarrow{\raisebox{-0.2ex}{\ensuremath{{\scriptscriptstyle #1}}}}} 

\newcommand{\boldemph}[1]{\emph{\textbf{#1}}} 

\newtheorem{theorem}{Theorem}[section]

\newtheorem*{theorem*}{Theorem}

\newtheorem{theoremalpha}{Theorem}

  \newaliascnt{proposition}{theorem}
  \newtheorem{proposition}[proposition]{Proposition}
  \aliascntresetthe{proposition}

  \newaliascnt{lemma}{theorem}
  \newtheorem{lemma}[lemma]{Lemma}
  \aliascntresetthe{lemma}

  \newaliascnt{corollary}{theorem}
  \newtheorem{corollary}[corollary]{Corollary}
  \aliascntresetthe{corollary}

  \newaliascnt{conjecture}{theorem}
  \newtheorem{conjecture}[conjecture]{Conjecture}
  \aliascntresetthe{conjecture}

\theoremstyle{definition}

  \newaliascnt{definition}{theorem}
  \newtheorem{definition}[definition]{Definition}
  \aliascntresetthe{definition}

  \newaliascnt{remark}{theorem}
  \newtheorem{remark}[remark]{Remark}
  \aliascntresetthe{remark}

  \newaliascnt{question}{theorem}
  \newtheorem{question}[question]{Question}
  \aliascntresetthe{question}

  \newaliascnt{example}{theorem}
  \newtheorem{example}[example]{Example}
  \aliascntresetthe{example}

\newcommand{\hyref}[1]{\hyperref[#1]{\ref{#1}}}

\newcommand{\step}[1]{\smallskip\noindent\emph{#1}}
\newcommand{\Case}[1]{\step{Case #1:}}
\newcommand{\Subcase}[1]{\step{Subcase #1:}}
\newcommand{\iffstep}[1]{\step{On (#1):}}
\newcommand{\Step}[1]{\step{Step #1:}}

\newcommand{\introitem}[2]{\smallskip\noindent\textbf{#1} (#2).}

%% file: header-tikz-surfaces.tex
\newcommand*{\tem}{0.8em}  
\newcommand*{\tpt}{0.8pt}  
\newcommand*{\trt}{1.125em}  

\tikzset{curve/.style={shape=circle,draw,fill=white,inner sep=\tpt,minimum width=\tem}}
\tikzset{ellipse/.style={shape=ellipse,draw,fill=white,inner sep=\tpt,minimum width=\tem,minimum height=0.85em}}
\tikzset{tripl/.style={shape=circle,draw,fill=white,minimum width=\trt}}
\tikzset{ldots/.style={shape=circle,fill=white,inner sep=\tpt,minimum width=\tem}}
\tikzset{shadow/.style={shape=ellipse,fill=white,minimum height=\trt,inner sep=0pt,text width=1pt}} 

\newcommand{\intext}[1]{\smash{\raisebox{-0.2ex}{#1}}}

\newcommand{\printdots}{\tiny\ldots}
\newcommand{\tinyminus}{\raisebox{0.2ex}{-}\hspace{-0.05em}}
\newcommand{\printnum}[1]{\tiny\sf \hspace{-0.075em}\tinyminus #1}

\newcommand{\tikzheader}[2]{\scalebox{1.08}{\ensuremath{\raisebox{-0.5ex}{\tikz[#1]{#2}}}}}

\pgfdeclarelayer{bg}    
\pgfsetlayers{bg,main}  

\newcommand\tikzchainsign[2]{
  \tikzheader{}{%
    \def\intnumlist{#1}
    \def\signlist{#2}
    \foreach [count=\x, evaluate=\x as \s using {{\signlist}[\x-1]}] \intnum in \intnumlist
    {
      \ifthenelse {\equal{\intnum}{}} 
             {\node[ldots] (p\x) at (0.6*\x-0.6,0) {\printdots};}
             {\ifthenelse {\equal{\s}{1}} %
                  {\node[curve] (p\x) at (0.6*\x-0.6,0) {\printnum{\intnum}};}
                  {\node[curve] (p\x) at (0.6*\x-0.6,0) {\tiny\sf \intnum};}}
    }
    \begin{pgfonlayer}{bg}
      \node[shadow] at (0,0) {}; 
      \draw (p1) to (p\x) ;
    \end{pgfonlayer}
  }
}

\newcommand\tikzchainpositive[1]{
  \tikzheader{}{%
    \def\intnumlist{#1}
    \foreach [count=\x] \intnum in \intnumlist
    {
      \ifthenelse {\equal{\intnum}{}} 
             {\node[ldots] (p\x) at (0.6*\x-0.6,0) {\printdots};}
             {\node[curve] (p\x) at (0.6*\x-0.6,0) {\tiny\sf \intnum};}
    }
    \begin{pgfonlayer}{bg}
      \node[shadow] at (0,0) {}; 
      \draw (p1) to (p\x) ;
    \end{pgfonlayer}
  }
}

\newcommand\tikzchain[1]{
  \tikzheader{}{%
    \def\intnumlist{#1}
    \foreach [count=\x] \intnum in \intnumlist
    {
      \ifthenelse {\equal{\intnum}{}} 
             {\node[ldots] (p\x) at (0.6*\x-0.6,0) {\printdots};}
             {\node[curve] (p\x) at (0.6*\x-0.6,0) {\printnum{\intnum}};}
    }
    \begin{pgfonlayer}{bg}
      \node[shadow] at (0,0) {}; 
      \draw (p1) to (p\x) ;
    \end{pgfonlayer}
  }
}

\newcommand\tikzsimplechainpure[2]{%
  \tikzheader{}{%
    \def\intnumlist{2,,2,1}
    \node[ellipse] (p0) at (-0.7,0) {\printnum{#1}#2};
    \foreach [count=\x] \intnum in \intnumlist
    {
      \ifthenelse {\equal{\intnum}{}} 
             {\node[ldots] (p\x) at (0.6*\x-0.6,0) {\printdots};}
             {\node[curve] (p\x) at (0.6*\x-0.6,0) {\printnum{\intnum}};}
    }
    \begin{pgfonlayer}{bg}
      \node[shadow] at (0,0) {}; 
      \draw (p0) to (p\x) ;
    \end{pgfonlayer}
  }
}

\newcommand\tikzsimplechain[1]{\tikzsimplechainpure{#1}{\tinyminus1}}

\newcommand\tikzchainnonred[2]{
  \tikzheader{ lbend/.style={bend left=10}, rbend/.style={bend right=10}}%
  {
    \def\intnumlist{#1}
    \def\multlist{#2}   
    \foreach [count=\x, evaluate=\x as \m using {{\multlist}[\x-1]}] \intnum in \intnumlist 
      {
      \node[curve] (p\x) at (0.6*\x-0.6,0) {\printnum{\intnum}};
      \begin{pgfonlayer}{bg}
      \node[shadow] at (p1) {};
        \ifthenelse {\equal{\m}{} \OR \equal{\m}{0}}
          {}
          {
            \ifthenelse {\equal{\m}{1}}
               { \draw (0.6*\x-0.6,0) to (0.6*\x,0) ;}         
               { \draw[style=lbend] (0.6*\x-0.6,0) to (0.6*\x,0) ; 
                 \draw[style=rbend] (0.6*\x-0.6,0) to (0.6*\x,0) ;}   
          }
      \end{pgfonlayer}
      }
  }
}

\newcommand\tikzchainmult[2]{
  \tikzheader{}{%
    \def\intnumlist{#1}
    \def\multlist{#2}
    \foreach [count=\x, evaluate=\x as \m using {{\multlist}[\x-1]}] \intnum in \intnumlist 
    {
      \ifthenelse {\equal{\intnum}{}}
        {\node[ldots] (p\x) at (0.6*\x-0.6,0) {\printdots};}
        {\ifthenelse {\equal{\m}{2}}
          {\node[curve,double] (p\x) at (0.6*\x-0.6,0) {\printnum{\intnum}};}
          {\ifthenelse {\equal{\m}{3}}
            {
              \node[tripl] at (0.6*\x-0.6,0) {};
              \node[curve,double] (p\x) at (0.6*\x-0.6,0) {\printnum{\intnum}};
            }
            {
              \node[curve] (p\x) at (0.6*\x-0.6,0) {\printnum{\intnum}};
            }
          }
        }
    }
    \begin{pgfonlayer}{bg}
      \node[shadow] at (0,0) {}; 
      \draw (p1) to (p\x);
    \end{pgfonlayer}
  }
}

\newcommand\tikztree[2]{
  \tikzheader{}{%
    \def\intnumlist{#1}
    \foreach [count=\x] \intnum in \intnumlist 
    {
      \ifthenelse {\equal{\intnum}{}}
           {\node[ldots] (p\x) at (0.6*\x-0.6,0) {\printdots};}
           {\node[curve] (p\x) at (0.6*\x-0.6,0) {\printnum{\intnum}};}
    }
    \node[curve] (q) at (0.6*\x/2-0.6*1/2,0.95*1/2) {\printnum{#2}};
    \begin{pgfonlayer}{bg}
      \foreach [count=\x] \intnum in {#1} 
      {
        \ifthenelse {\equal{\intnum}{}} 
           {} 
           {\draw (q) -- (p\x);}
      }
    \end{pgfonlayer}
  }
}

\newcommand\tikztreecurved[2]{
  \tikzheader{}{%
    \def\intnumlist{#1}
    \foreach [count=\x] \intnum in \intnumlist
    {
      \ifthenelse {\equal{\intnum}{}}
           {\node[ldots] (p\x) at (0.6*\x-0.6,0) {\printdots};}
           {\node[curve] (p\x) at (0.6*\x-0.6,0) {\printnum{\intnum}};}
    }
    \node[curve] (q) at (-0.6,0.1) {\printnum{#2}};
    \begin{pgfonlayer}{bg}
      \foreach [count=\x] \intnum in {#1}
      {
        \ifthenelse {\equal{\intnum}{}}
           {}
           {\draw ($(q)+(0,-0.05+\x/20)$) to [out=4*\x,in=150] ($(p\x)+(0,\x/30)$);}
      }
    \end{pgfonlayer}
  }
}

\newcommand\tikztreemult[4]{
  \tikzheader{}{%
    \def\intnumlist{#1}
    \def\multlist{#3}
    \foreach [count=\x, evaluate=\x as \m using {{\multlist}[\x-1]}] \intnum in \intnumlist
    {
      \ifthenelse {\equal{\intnum}{}} 
        {\node[ldots] (p\x) at (0.6*\x-0.6,0) {\printdots};} 
        {\ifthenelse {\equal{\m}{2}}
          { \node[curve,double] (p\x) at (0.6*\x-0.6,0) {\printnum{\intnum}}; }
          { \ifthenelse {\equal{\m}{3} }
              {
                \node[tripl] at (0.6*\x-0.6,0) {};
                \node[curve,double] (p\x) at (0.6*\x-0.6,0) {\printnum{\intnum}};
              }
              {
                \node[curve] (p\x) at (0.6*\x-0.6,0) {\printnum{\intnum}};
              }
          }
        }
    }
    \ifthenelse {\equal{#4}{2}} 
      {\node[curve,double] (q) at (0.6*\x/2-0.6*1/2,0.95*1/2) {\printnum{#2}};}
      {\ifthenelse {\equal{#4}{3}}
        {
          \node[tripl] at (0.6*\x/2-0.6*1/2,0.85*1/2) {};
          \node[curve,double] (q) at (0.6*\x/2-0.6*1/2,0.95*1/2) {\printnum{#2}};
        }
        {
          \node[curve] (q) at (0.6*\x/2-0.6*1/2,0.75*1/2) {\printnum{#2}};
        }
      }
    \begin{pgfonlayer}{bg}
      \foreach [count=\x] \intnum in {#1} 
      {
        \ifthenelse {\equal{\intnum}{}} 
          {}
          {\draw (q) -- (p\x);}
      }
    \end{pgfonlayer}
  }
}

\newcommand\tikztypeD[4]{
  \tikzheader{}{%
    \def\intnumlist{#3}
    \def\multlist{#4}
    \node[curve] (upperleaf) at (-0.5,-0.3) {\printnum{#1}};
    \node[curve] (lowerleaf) at (-0.5, 0.3) {\printnum{#2}};
    \foreach [count=\x, evaluate=\x as \m using {{\multlist}[\x-1]}] \intnum in \intnumlist
    {
      \ifthenelse {\equal{\intnum}{}} 
        {\node[ldots] (p\x) at (0.6*\x-0.6,0) {\printdots};} 
        {\ifthenelse {\equal{\m}{2}}
          { \node[curve,double] (p\x) at (0.6*\x-0.6,0) {\printnum{\intnum}}; }
          { \ifthenelse {\equal{\m}{3} }
              {
                \node[tripl] at (0.6*\x-0.6,0) {};
                \node[curve,double] (p\x) at (0.6*\x-0.6,0) {\printnum{\intnum}};
              }
              {
                \node[curve] (p\x) at (0.6*\x-0.6,0) {\printnum{\intnum}};
            }
          }
        }
      }
    \begin{pgfonlayer}{bg}
      \draw (p1) -- (p\x);
      \draw (upperleaf) -- (p1) -- (lowerleaf);
    \end{pgfonlayer}
  }
}

\newcommand\tikztypeDmirror[4]{
  \tikzheader{}{%
    \def\intnumlist{#3}
    \def\multlist{#4}
    \node[curve] (upperleaf) at (0.5,-0.3) {\printnum{#1}};
    \node[curve] (lowerleaf) at (0.5, 0.3) {\printnum{#2}};
    \foreach [count=\x, evaluate=\x as \m using {{\multlist}[\x-1]}] \intnum in \intnumlist
    {
      \ifthenelse {\equal{\intnum}{}} 
        {\node[ldots] (p\x) at (-0.6*\x+0.6,0) {\printdots};} 
        {\ifthenelse {\equal{\m}{2}}
          { \node[curve,double] (p\x) at (-0.6*\x+0.6,0) {\printnum{\intnum}}; }
          { \ifthenelse {\equal{\m}{3} }
              {
                \node[tripl] at (-0.6*\x+0.6,0) {};
                \node[curve,double] (p\x) at (-0.6*\x+0.6,0) {\printnum{\intnum}};
              }
              {
                \node[curve] (p\x) at (-0.6*\x+0.6,0) {\printnum{\intnum}};
              }
          }
        }
    }
    \begin{pgfonlayer}{bg}
      \draw (p1) -- (p\x);
      \draw (upperleaf) -- (p1) -- (lowerleaf);
    \end{pgfonlayer}
  }
}

%% file: biblio-surfaces.tex
\newcommand{\arXiv}[1]{{\href{http://arxiv.org/abs/#1}{\texttt{arXiv:#1}}}}
\newcommand{\jstor}[1]{{\href{http://www.jstor.org/stable/#1}{JSTOR}}}
\newcommand{\eudml}[1]{{\href{https://eudml.org/doc/#1}{EuDML}}}

\begin{spacing}{0.9}

\end{spacing}